\documentclass[preprint,sort&compress,final,12pt]{elsarticle}

\usepackage{amsmath}
\usepackage{amssymb}

\newtheorem{theorem}{Theorem}[section]

\newtheorem{definition}[theorem]{Definition}
\newtheorem{lemma}[theorem]{Lemma}
\newtheorem{proposition}[theorem]{Proposition}

\numberwithin{equation}{section}

\def\Proof{\noindent{\bf Proof.}~}
\def\qed{\hfill$\square$\smallskip}
\def\dint{\displaystyle\int}

\def\dint{\displaystyle\int}

\def\dfrac#1#2{\frac{\displaystyle {#1}}{\displaystyle {#2}}}
\def\Re{\mathrm{Re}}

\def\Res{\mathrm{Res}}

\def\R{{\Bbb R}}
\def\C{{\Bbb C}}
\journal{Stochastic Processes and their Applications}

\begin{document}

\begin{frontmatter}

\title{Second Moment Boundedness of Linear Stochastic Delay Differential Equations}

\author[au1]{Zhen Wang}

\author[au1]{Xiong Li\footnote{Supported by the Natural Science Foundation of China (NSFC 11031002), the Fundamental
Research Funds for the Central Universities and the Scientific Research Foundation
for the Returned Overseas Chinese Scholars, State Education Ministry.}}

\address[au1]{School of Mathematical Sciences, Beijing Normal University, Beijing 100875, P.R. China.}

\author[au3]{Jinzhi Lei\footnote{Supported by the Natural Science Foundation of China  (NSFC 11272169).}}
\address[au3]{Zhou Pei-Yuan Center for Applied Mathematics, Tsinghua University, Beijing 100084, P.R. China.}
%\ead[au3]{jzlei@mail.tsinghua.edu.cn}

\begin{abstract}
This paper studies the second moment boundedness of solutions of linear stochastic delay differential equations. First, we give a framework, for general $\mathrm{N}$-dimensional linear stochastic differential equations with a single discrete delay, of calculating the characteristic function for the second moment boundedness. Next, we apply the proposed framework to a special case of a type of $2$-dimensional equation that the
stochastic terms are decoupled. For the $2$-dimensional equation, we obtain the characteristic function explicitly given by equation coefficients, the characteristic function gives sufficient conditions for the second moment to be bounded or unbounded.
\end{abstract}


\end{frontmatter}

\section{Introduction}

Stochastic delay differential equations have been extensively studied in the last several decades from different points of view
(see \cite{Caraballo10,Duan04,Ito,Ivanov03,LeiM,Mackey94,Mackey95,Mao1,Mao2,Mao3,Mohammed03,Mohammed04} and the references therein).
However, many basic issues remain unsolved even for linear equations with constant coefficients.

In this paper, we study general $N$-dimensional linear stochastic delay differential equations with a single discrete delay
(here we assume the delay $\tau=1$)
\begin{equation}
\label{sdde2}
{\small
\begin{array}{c}
dx_i(t) = \left[a_i^j x_j(t) + b_i^j x_j(t-1)\right] d t + \left[\mu_i^k + \sigma_i^{jk}x_j(t) + \eta_i^{jk}x_j(t-1)\right]dW_{k}.\\
(i, j=1,\cdots, N; \quad k=1,\cdots, K)
\end{array}}
\end{equation}
Here the \textit{Einstein summation convention} has been used so that repeated indices are implicitly summed over, all coefficients $a_i^j$, $b_i^j$,
$\mu_i^k$, $\sigma_i^{jk}$, $\eta_i^{jk}$ are constants, and $W_k$ are independent $1$-dimensional Wiener Processes. The initial functions
are assumed to be $x_i=\phi_i\in C([-1,0],\R)\,(i=1,\cdots, N)$. This paper considers the second moment boundedness of solutions of \eqref{sdde2}.
We always assume It\^{o} interpretation for the stochastic integral, and results for Stratonovich interpretation can be
obtained similarly.

To the best of our knowledge, there are very few results for the stability and second moment boundedness for equation \eqref{sdde2}, and
most of known results are obtained through the method of Lyapunov functional \cite{Mackey94,Mackey95,Mao1,Mohammed84}.
However, it remains unclear how can we define the characteristic equation  for the boundedness of
the solution moments of \eqref{sdde2}. In 2007, Lei and Mackey \cite{LeiM} introduced the method of Laplace transform to study the second moment
boundedness of $1$-dimensional equations with a single discrete delay, and proposed a
characteristic equation. In \cite{WLL}, the authors extended the method to the $1$-dimensional equation with distributed delay.

Based on previous studies \cite{LeiM,WLL}, this paper aims at proposing  a general framework to calculate the characteristic function for high dimensional
linear stochastic delay differential equations with a single delay. The obtained characteristic function (as we can see below) is complicate, that shows the elaborate correlations when delay and stochastic effects are coupled into an dynamical system. As an example,
we apply the framework to study a special situation of a $2$-dimensional equation in which the stochastic terms are decoupled
($N=K=2$, and $\mu_i^k = \sigma_i^{jk} = \eta_i^{jk} = 0$ when $i\not= k$).

Rest of this paper is organized as follows. In Section 2, we briefly introduce basic results for linear delay differential equations. In Section 3,
a general framework for defining the characteristic function of a $\mathrm{N}$-dimensional stochastic delay differential equation is given.
In Section 4, we discuss the boundedness of the second moment of \eqref{sdde2} for a simple case: $N=K=2$ and the stochastic terms are decoupled.
Theorem \ref{Thm-unb} establishes the unbounded condition for the second moment if the trivial solution of the unperturbed equation
is unstable. When the trivial solution of the unperturbed equation is stable, we obtain a characteristic function given explicitly by the equation
parameters. Boundedness of the second moments depends on the supremum of the real parts of all roots of the characteristic equation
(Theorem \ref{Thm-bound}). An explicit condition for the boundedness of the second moment is also proved following framework of calculation given here (Theorem \ref{Thm-BK}).

\section{Preliminaries}

In this paper, we always use the $L^1$ norm for a tensor. For example, the $L^1$ norm of a second order tensor $A = \{a_i^j\}$ is
\begin{equation}
\|A\| = \sum_{i,j=1}^{N}|a_i^j|.
\end{equation}
For $\phi=(\phi_i)_{N \times 1} \in C([-1,0],\R^N) $, the norm is defined as
\begin{equation}
\|\phi\| = \sup_{\theta\in [-1,0]}\sum_{i=1}^{N} |\phi_i(\theta)|.
\end{equation}

In this section we give some basic results for the $N$-dimensional linear delay differential
equation
\begin{equation}\label{dde}
\dfrac{d x_i(t)}{dt}= a_i^j x_j(t) + b_i^j x_j(t-1) \;\;(i,j=1,\cdots, N )
\end{equation}
with initial functions $x_i=\phi_i\in C([-1,0],\R)$. The linear autonomous functional differential equation \eqref{dde} has been studied
extensively and details can be referred to \cite{Arino,Bellman,Hale93}.

The fundamental matrix of \eqref{dde}, denoted by
$$
X(t)=\left(X_i^j(t)\right)_{N\times N},
$$
is the solution of \eqref{dde} with the initial condition (hereinafter $\delta$ means the Kronecker delta)
$$
X_i^j(t)=\left\{\begin{array}{cc}
\delta_i^j, & t=0, \\ 0, &  \quad -1\leq t <0.
\end{array}\right.
$$
Using the fundamental matrix $X(t)$, the solution of \eqref{dde} with initial function $\phi \in C([-1,0],\R^N)$ can be represented as
\begin{equation}
\label{xphi}
x_{\phi,i}(t) = X_i^j(t)\phi_j(0)+\int^{0}_{-1}X_i^j(t-1-\theta)b_j^l\phi_l(\theta)d\theta.
\end{equation}
From \eqref{xphi}, the asymptotic behavior of $x_{\phi}(t)=\left(x_{\phi,i}(t)\right)_{N\times 1}$ is determined by the fundamental matrix $X(t)$.

Denote the matrices
\begin{equation}\label{matrices}
A=(a_i^j)_{N\times N},\quad  B=(b_i^j)_{N\times N}, \quad  I =(\delta_i^j)_{N\times N}, \quad  \mu =(\mu_i^j)_{N\times N}.
\end{equation}
Taking the Laplace transform on both sides of \eqref{dde}, we obtain
\begin{equation}
\label{LY}
\mathcal{L}(X)(\lambda) = [\Delta(\lambda)]^{-1},
\end{equation}
where
\begin{equation}
\label{Delta}
\Delta(\lambda) = \lambda I - A - B e^{-\lambda}.
\end{equation}
Thus, $h(\lambda) = \det(\Delta(\lambda))$ is the characteristic function for the linear stability of \eqref{dde}.

The following results are straightforward from the above discussions.

\begin{theorem} \label{Thm-X}
Let
\begin{equation}\label{alpha0}
\alpha_{0}=\sup\{\Re(\lambda):  h(\lambda)=0,\, \lambda \in\C\}.
\end{equation}
Then
\begin{enumerate}
\item[\textnormal{(i)}] for any $\alpha>\alpha_{0}$ there exists a constant $\bar{K}=\bar{K}(\alpha)\geq 1$ such that
the fundamental matrix $X(t)$ satisfies
\begin{equation*}
\|X(t)\|\leq \bar{K} e^{\alpha t}, \quad t\geq 0;
\end{equation*}
\item[\textnormal{(ii)}] for any $\alpha>\alpha_{0}$ there exists a constant
$\tilde{K}=\tilde{K}(\alpha)\geq 1$ such that for any $\phi\in C([-1,0],\R^N)$
the solution $x_{\phi}(t)$ of (\ref{dde}) satisfies
\begin{equation*}
\|x_{\phi}(t)\|\leq \tilde{K}\|\phi\|e^{\alpha t}, \quad  t\geq 0;
\end{equation*}
\item[\textnormal{(iii)}] for any $\alpha_{1}<\alpha_{0}$, there exists $\bar{\alpha}\in (\alpha_{1}, \alpha_{0})$ and a
subset $U \subset \R^{+} $ with measure $m(U)=+\infty$ such that for any $i,j=1,\cdots, N$,
\begin{equation}\label{Xij}
\|X_{i}^{j}(t)\| \geq e^{\bar{\alpha}t},\quad t\in U.
\end{equation}
\end{enumerate}
\end{theorem}

Here, the number $\alpha_0$ is also termed as the \textit{Lyapunov exponent} (Definition 1.19 in \cite{Arino}). The proofs of (i) and (ii) in
Theorem \ref{Thm-X} is referred to the proofs in \cite[Theorem 1.21 in Chapter 3]{Arino}. The proof of (iii) is similar to that of
\cite[Theorem 2.3, 2]{WLL} and is detailed at Appendix A.

From Theorem \ref{Thm-X}, the trivial solution of \eqref{dde} is locally asymptotically stable if and only if $\alpha_{0}<0.$

\section{Moments of linear stochastic delay differential equations--General cases}

Now, we discuss the solution moments and a framework for calculating the  characteristic function of the second moment of solutions for general cases.

The existence and uniqueness results for stochastic delay differential  equations have been established in \cite{Ito,Mao1,Mohammed84}.
Using the fundamental matrix $X(t)$,  the solution of \eqref{sdde2} with initial function
$x=\phi\in C([-1,0],\R^N)$ is a $N$-dimensional stochastic process given by It\^{o}
integral as follows
\begin{equation}\label{xtphi}
x_i(t;\phi)=x_{\phi,i}(t) + \int^{t}_{0} X_i^l(t-s)\big(\mu_l^k + \sigma_l^{jk}x_j(s;\phi) + \eta_l^{jk}x_j(s-1;\phi)\big) dW_k,\; t\geq 0,
\end{equation}
where $(x_{\phi,i}(t))_{N \times 1}$ is the solution of \eqref{dde} defined by \eqref{xphi}.

We denote by $E$ the mathematical expectation. Now we give the definitions of the $p^{\mathrm{th}}$ moment exponential stability and the $p^{\mathrm{th}}$ moment boundedness.

\begin{definition}\label{def3.1}
A solution $x(t;\phi)$ of \eqref{sdde2} is said to be the first moment exponentially
stable if there exist two positive constants $\gamma$ and $R$ such that
\begin{equation*}
\|E(x(t;\phi))\|\leq R\|\phi\|e^{-\gamma t},\quad t\geq 0,
\end{equation*}
for all $\phi\in C((-1,0],\R^N)$. When $p\geq 2$, a solution of
\eqref{sdde2} is said to be the $p^{\mathrm{th}}$ moment exponentially stable if there
exist two positive constants $\gamma$ and $R$ such that
\begin{equation*}
E\big(\|x(t;\phi)-E(x(t;\phi))\|^{p}\big)\leq R\|\phi\|^{p}e^{-\gamma t},\quad t\geq 0,
\end{equation*}
for all $\phi\in C((-1,0],\R^N)$.
\end{definition}

\begin{definition}\label{def3.2}
For $p\geq 2$, a solution $x(t;\phi)$ of \eqref{sdde2} is said to be the $p^{\mathrm{th}}$ moment
bounded if there exists a positive constant $\tilde{C}=\tilde{C}(\|\phi\|^p)$
such that
\begin{equation*}
E\big(\|x(t;\phi)-E(x_i(t;\phi))\|^{p}\big)\leq \tilde{C},\quad t\geq 0,
\end{equation*}
for all $\phi\in C((-1,0],\R^N)$. Otherwise, the $p^{\mathrm{th}}$ moment
is said to be unbounded.
\end{definition}

\subsection{First moment}

From \eqref{xtphi} and applying the properties of It\^{o} integral, it is easy to have
\begin{equation*}
E x_i(t; \phi)=x_{\phi,i}(t)\; (i=1, \cdots, N).
\end{equation*}
Thus, the following result is straightforward from Theorem \ref{Thm-X}.

\begin{theorem}
\label{thm:fm}
Let $\alpha_0$ be defined by \eqref{alpha0}. For any $\alpha > \alpha_0$, there exists a positive constant $\tilde{K}$
(defined as Theorem \ref{Thm-X}) such that for any $\phi\in C([-1,0],\mathbb{R}^N)$, the solution
$x(t; \phi)=(x_{i}(t;\phi))_{N \times 1}$ of \eqref{sdde2} satisfies
\begin{equation}\label{Ex}
\|E x(t; \phi)\| \leq \tilde{K}\|\phi\|e^{\alpha t},\quad t\geq 0.
\end{equation}
In particular, the first moment of \eqref{sdde2} is exponentially stable when $\alpha_{0}<0$.
\end{theorem}

\subsection{Second moment}

Now we study the second moment. First, we give some notations. Let $x(t;\phi) \triangleq (x_i(t))_{N \times 1}$
be a solution of \eqref{sdde2} and $ \tilde{x}_i(t) = x_i(t) - E(x_i(t))\; (i=1, \cdots, N)$, and define
\begin{equation}
M_{ij}(t) = E(\tilde{x}_i(t)\tilde{x}_j(t)),\quad N_{ij}(t) = E(\tilde{x}_i(t) \tilde{x}_j(t-1)).
\end{equation}
Then $M_{ii}(t)$ is the second moment of $x_i(t)$. It is easy to have $\tilde{x}_i(t)=M_{ij}(t)=N_{ij}(t)=0$ when $t \in[-1, 0]$, and $E(\tilde{x}_i(t))=0$ for all $t\geq 0$.

Denote
\begin{equation}
\Sigma_l^k(s) = \mu_l^k + \sigma_l^{jk}x_j(s) + \eta_l^{jk}x_j(s-1).
\end{equation}
From \eqref{xtphi}, we have
$$
\tilde{x}_i(t) = \int_0^t X_i^l(t-s) \Sigma_l^k(s)d W_k.
$$
Since $E(dW_k dW_m) = \delta_{k,m}$, we obtain
\begin{equation}
M_{ij}(t) = E(\tilde{x}_i(t)\tilde{x}_j(t)) = \int_0^t X_i^l(t-s)E\left[\Sigma_l^k(s) \delta_{k,m} \Sigma_{p}^m(s)\right] X_j^p(t-s) d s.
\end{equation}

\subsubsection{Additive noise}

We have additive noise when $\sigma_l^{jk} = \eta_l^{jk} = 0$ for any $i,j,k$  and $\mu \not=0$ ($\mu$ is defined at \eqref{matrices}). In this case, we have
\begin{equation}\label{Mijt2}
M_{ij}(t) = \int_0^t X_i^l(s)X_j^p(s)\mu_l^k\delta_{k,m}\mu_p^m ds.
\end{equation}
Thus the upper bound of $M(t)$ is determined by that of the fundamental matrix $X(t)$. Hence, we have the following sufficient conditions for the second moments $M(t)$ to be bounded or unbounded.

\begin{theorem}\label{Thm-bou}
Let $\alpha_{0}$ be defined as \eqref{alpha0}. Assume $\sigma_l^{jk} = \eta_l^{jk} = 0$ for any $i,j,k$ and $\mu \not=0$.
Then the second moment of  \eqref{sdde2} is bounded if $\alpha_{0}<0$, and unbounded if $\alpha_{0}>0$.
\end{theorem}

The proof is similar to that in \cite[Theorem 3.4]{WLL} and is omitted here.

The critical case $\alpha_{0}=0$ is not discussed here and the boundedness issue remains open.

\subsubsection{A general framework}

Now, we consider the general situation. Let
\begin{eqnarray*}
P_{lp}(s) &=& E(\Sigma_l^k(s))\delta_{k,m}E(\Sigma_p^m(s)), \\
Q_{lp}(s) &=& E\left[ \left(\sigma_l^{jk}\tilde{x}_j(s) + \eta_l^{jk}\tilde{x}_j(s-1) \right) \delta_{k,m}
\big(\sigma_p^{qm}\tilde{x}_q(s) + \eta_p^{qm}\tilde{x}_q(s-1) \big) \right].
\end{eqnarray*}
Then $P_{lp}(s)= P_{pl}(s),\; Q_{lp}(s)= Q_{pl}(s)$ and
\begin{eqnarray} \label{Qlp}
Q_{lp}(s)
&=& \sigma_l^{jk}\delta_{k,m}\sigma_p^{qm} M_{jq}(s) + \eta_l^{jk} \delta_{k,m}\eta_p^{qm} M_{jq}(s-1) \nonumber \\
&&{} + \left(\sigma_l^{jk}\delta_{k,m}\eta_p^{qm} + \eta_l^{qk}\delta_{k,m}\sigma_p^{jm} \right) N_{jq}(s).
\end{eqnarray}
We note that $x_i(s)=\tilde{x}_i(s) + E x_i(s)$, then
$$
E\left[\Sigma_l^k(s)\delta_{k,m}\Sigma_p^m(s)\right] = P_{lp}(s) + Q_{lp}(s).
$$
Therefore,
$$
M_{ij}(t) = \int_0^t X_i^l(t-s)X_j^p(t-s)(P_{lp}(s) + Q_{lp}(s))d s.
$$
Denote
\begin{equation}\label{Ft}
F_{ij}(t) = \int_0^t X_i^l(t-s)X_j^p(t-s) P_{lp}(s)d s,
\end{equation}
then
\begin{equation}\label{Mt}
M_{ij}(t) = F_{ij}(t) + \int_0^t X_i^l(t-s)X_j^p(t-s) Q_{lp}(s) d s.
\end{equation}
Similarly, we have
\begin{equation}\label{Nt}
N_{ij}(t) = \int_0^{t-1} X_i^l(t-s)X_j^p(t-1-s)(P_{lp}(s)+Q_{lp}(s)) d s.
\end{equation}

From \eqref{Mt}, we have $|M_{ii}(t)| \geq |F_{ii}(t)|\ (\forall i)$. If $\alpha_0 > 0$ ($\alpha_0$ is defined by \eqref{alpha0}), similar to discussions in \cite{WLL}, $|F_{ii}(t)|$ approaches to infinity exponentially, and therefore the second moment is unbounded (a proof for the case of a $2$-dimensional equation is given in the next section). Thus, we only need to study the situation when $\alpha_0 < 0$.

From \eqref{Mt} and \eqref{Nt}, and take Laplacians to both sides of them (existence of the Laplacians are proved in Lemmas \ref{Lem-FL}
and \ref{Lem-ML} below), we obtain
\begin{equation}
\label{LM}
\mathcal{L}(M_{ij}) = \mathcal{L}(X_i^hX_j^p)\left[\mathcal{L}(P_{hp}) + \mathcal{L}(Q_{hp})\right],
\end{equation}
and
\begin{equation}
\label{LN}
\mathcal{L}(N_{ij}) = \mathcal{L}(X_i^h(t)X_j^p(t-1))\left[\mathcal{L}(P_{hp}) + \mathcal{L}(Q_{hp})\right].
\end{equation}
From \eqref{Qlp}, we have
\begin{eqnarray}
\label{LQ}
\mathcal{L}(Q_{hp}) &=& (\sigma_h^{jk}\delta_{k,m}\sigma_p^{qm} + e^{-\lambda}\eta_h^{jk}
\delta_{k,m}\eta_p^{qm}) \mathcal{L}(M_{jq}) \nonumber\\
&&{} + (\sigma_h^{jk}\delta_{k,m}\eta_p^{qm} + \eta_h^{qk}\delta_{k,m}\sigma_p^{jm})\mathcal{L}(N_{jq}).
\end{eqnarray}

Now, one can solve $\mathcal{L}(M_{ij})$ from \eqref{LM}-\eqref{LQ} following the procedure below. First, solve $\mathcal{L}(P_{hp}) + \mathcal{L}(Q_{hp})$ by $\mathcal{L}(M_{ij})$ from \eqref{LM} as
$$\mathcal{L}(P_{hp}) + \mathcal{L}(Q_{hp}) = S_{hp}^{ij}\mathcal{L}(M_{ij}).$$
Here $S_{hp}^{ij}$ is the inverse tensor of $\mathcal{L}(X_i^h X_j^p)$. Next, substitute the obtained $\mathcal{L}(P_{hp}) + \mathcal{L}(Q_{hp})$ into \eqref{LN} to linearly express $\mathcal{L}(N_{ij})$ through $\mathcal{L}(M_{ij})$:
$$\mathcal{L}(N_{ij})=\mathcal{L}(X_i^h(t)X_j^p(t-1))S_{hp}^{kl}\mathcal{L}(M_{kl}).$$
Then, put the resulting $\mathcal{L}(N_{ij})$ into \eqref{LQ} so that $\mathcal{L}(Q_{hp})$ linearly depends on $\mathcal{L}(M_{ij})$ in the form
$$\mathcal{L}(Q_{hp}) = T_{hp}^{kl}\mathcal{L}(M_{kl}).$$
Finally, substituting $\mathcal{L}(Q_{hp})$ back to \eqref{LM} to obtain an equation for $\mathcal{L}(M_{ij})$ of form
\begin{equation}
\label{eq:LM}
\left(\delta_{ij}^{kl} - \mathcal{L}(X_i^hX_j^p)T_{hp}^{kl}\right)\mathcal{L}(M_{kl}) = \mathcal{L}(X_i^kX_j^l)\mathcal{L}(P_{kl}),
\end{equation}
Then equation \eqref{eq:LM} is a linear equation of $\mathcal{L}(M_{kl})$. Thus, the determinant of the coefficients, denoted by
\begin{equation}\label{H}
H(\lambda) = \det\left[ (\delta_{ij}^{kl} - \mathcal{L}(X_i^hX_j^p)T_{hp}^{kl})\right],
\end{equation}
is the desired characteristic function.

We note that \eqref{eq:LM} contains $N^2$ linear equations. Nevertheless, we can simplify the calculation due to symmetry. For example, since
$\mathcal{L}(M_{kl})=\mathcal{L}(M_{lk})$ and $\mathcal{L}(P_{kl})=\mathcal{L}(P_{lk})$ in \eqref{eq:LM}, we only need to solve equations for
$\mathcal{L}(M_{kl})$ with $k\leq l$, and therefore have $N(N+1)/2$ equations.

The above procedure gives a general framework to obtain the characteristic function. However, it is too complicate to obtain an explicit expression
for general cases. In the next section, we study a $2$-dimensional equation with a specific form.

Denote the matrices
$$
M=(M_i^j(t))_{N\times N},\;\;  N=(N_i^j(t))_{N\times N}, \;\;  F =(F_i^j(t))_{N\times N}, \;\; Q =(Q_i^j(t))_{N\times N}.
$$
Before introducing the results for the $2$-dimensional equation, we give some estimates for $F(t)$, $M(t)$ and $N(t)$ for general situation.
These estimations ensure the existence of Laplace transforms of $F(t), M(t)$ and $N(t)$.
\begin{lemma}\label{Lem-FL}
Let $\alpha_{0}$ be defined at \eqref{alpha0} and assume $\alpha_0 < 0$. Then for any $\alpha\in(\alpha_{0},0),$
there exists a positive constant $K_{1}=K_{1}(\alpha,\phi)\,(\phi\in C([-1,0],\R^{N}))$
such that
\begin{equation}
\|F(t)\|\leq K_{1}(1-e^{2\alpha t}),\quad t>0. \label{FK1}
\end{equation}
\end{lemma}

\Proof
From Theorem \ref{thm:fm}, for any $\alpha \in (\alpha_0, 0)$, there exists a positive constant $\tilde{K} = \tilde{K}(\alpha)$ such that
\begin{eqnarray*}
E(\Sigma_l^k(s)) &<& \|\mu\| + \tilde{K} \sum_{j=1}^2 \left(|\sigma_l^{jk}| + e^{-\alpha}|\eta_l^{jk}| \right)
\|\phi\|e^{\alpha s}\\
&<& \|\mu\| + \tilde{K} (\|\sigma\| + e^{-\alpha}\|\eta\|) \|\phi\|e^{\alpha s}
\end{eqnarray*}
for any $k,l$. Here $\|\cdot\|$, as we mentioned before, mean the $L^1$ norm of a tensor. Hence,
$$
|P_{lp}(s)|\leq \left(\|\mu\| + \tilde{K} (\|\sigma\| + e^{-\alpha}\|\eta\|) \|\phi\|e^{\alpha s}\right)^2.
$$
From Theorem \ref{Thm-X}, there exists a positive constant $\bar{K} = \bar{K(}\alpha)$ such that $\|X(t)\|< \bar{K} e^{\alpha t}$. Therefore, from \eqref{Ft}, it is not difficult to have a constant $K_1$ (defined by $\tilde{K}$ and $\bar{K}$) so that \eqref{FK1} holds.

\begin{lemma}\label{Lem-ML}
Let $\alpha_{0}$ be defined at \eqref{alpha0} and assume $\alpha_0<0$. For any $\alpha\in(\alpha_{0},0)$, there exists $\nu= \nu(\alpha)$ such that
\begin{equation}\label{MK1}
\|M(t)\| \leq K_{1} e^{\nu t}, \quad t\geq 0,
\end{equation}
and
\begin{equation}\label{NK1}
\|N(t)\| \leq \left(1+e^{-\nu}\right) K_1
e^{\nu t},\quad t\geq 0,
\end{equation}
where $K_{1}$ is defined as in Lemma \ref{Lem-FL}.
\end{lemma}

\Proof
From the Cauchy-Schwarz inequality, we obtain for $i,j=1,  \cdots, N$,
\begin{equation} \label{Nij}
|N_{ij}(t)|=|E\big(\tilde{x}_{i}(t)\tilde{x}_{j}(t-1)\big)| \leq \dfrac{M_{ii}(t)+M_{jj}(t-1)}{2}.
\end{equation}
Then
\begin{equation}
\label{NM}
\|N(t)\|  \leq \dfrac{1}{2}\sum_{i,j=1}^{N} \left (M_{ii}(t)+M_{jj}(t-1) \right) \leq \|M(t)\| + \|M(t-1)\|.
\end{equation}
Hence, from \eqref{Qlp} and \eqref{NM}, there exists $C_0 >0$ so that
\begin{eqnarray}\label{QNorm}
|Q_{lp}(s)| & \leq  & \left|\sigma_l^{jk} \delta_{k,m}\sigma_p^{qm} M_{jq}(s)\right|
+ \left|\eta_l^{jk}\delta_{k,m}\sigma_p^{qm}M_{jq}(s-1)\right| \nonumber\\
&&{} + \left|(\sigma_l^{jk}\delta_{k,m}\eta_p^{qm} + \eta_l^{qk}\delta_{k,m}\sigma_p^{jm}) N_{jq}(t)\right| \nonumber\\
&\leq & C_0 (\|M(s)\| + \|M(s-1)\|).
\end{eqnarray}
Thus, from Lemma \ref{Lem-FL}, we obtain for any $\alpha\in(\alpha_{0},0)$,
\begin{eqnarray*}
\|M(t)\| & \leq  & \|F(t)\| + \int_{0}^{t}\|X(t-s)\|^2 \|Q(s)\| d s  \\
& \leq &  K_{1}(1-e^{2\alpha t}) + C_0 \bar{K}^2 \int_{0}^{t} e^{2\alpha (t-s)} \left(\| M(s)\| + \| M(s-1)\|\right) d s\\
& \leq & K_{1} + C_0\bar{ K}^2 \int_{0}^{t} e^{2\alpha (t-s)}\| M(s)\| d s
+ C_0 \bar{K}^2 e^{-2 \alpha} \int_{-1}^{t-1} e^{2\alpha (t-s)} \| M(s)\| d s \\
& \leq & K_{1} + C_0 \bar{K}^2 ( 1 + e^{-2 \alpha} ) \int_{0}^{t} \| M(s)\| d s.
\end{eqnarray*}
Applying the Gronwall inequality, we have
\begin{equation*}
\|M(t)\| \leq K_{1} e^{\nu t},
\end{equation*}
where $\nu = C_0 \bar{K}^2 (1+e^{-2\alpha})$. The estimation \eqref{NK1} is obtained from \eqref{MK1} and \eqref{NM}.

\section{Application to $2$-dimensional equations}

In this section, we apply the general framework established in the above to a special case of a $2$-dimensional equation that $N=K=2$ and
$\mu_i^k = \sigma_i^{jk} = \eta_i^{jk} = 0$ when $i\not= k$. Hereafter, we do not use the Einstein summation convention, and
introduce following notations for simplicity: $\mu_i = \mu_i^i,\ \sigma_i^j=\sigma_{i}^{ji},\ \eta_i^j = \eta_i^{ji} (i,j=1,2)$.
Thus, the equation we studied becomes
\begin{equation} \label{sdde3}
\left\{\begin{array}{rcl}
dx_1(t) &=& \displaystyle\sum_{j=1}^2 \left(a_1^j x_j(t) + b_1^j x_j(t-1) \right) d t\\
&& \displaystyle{} + \left(\mu_1 +\sum_{j=1}^2 \left(\sigma_1^{j}x_j(t) + \eta_1^{j}x_j(t-1)\right) \right) dW_1, \\
dx_2(t) &=&\displaystyle\sum_{j=1}^2 \left(a_2^j x_j(t)+ b_2^j x_j(t-1)\right) d t\\
&& \displaystyle{} + \left(\mu_2  + \sum_{j=1}^2 \left(\sigma_2^{j}x_j(t) + \eta_2^{j}x_j(t-1)\right) \right) dW_{2},
\end{array}\right.
\end{equation}
In this particular case, the expressions of $P_{ij}(t)$, $Q_{ij}(t)$ and $F_{ij}(t)$, $M_{ij}(t)$, $N_{ij}(t)$ $(i,j= 1,2)$ in the previous section are as follows:
\begin{eqnarray*}
P_{12}(t) &=& P_{21}(t)=Q_{12}(t)=Q_{21}(t)=0,\\
P_{ii}(t) &=& \left(\mu_i + \sum_{j=1}^2\left(\sigma_i^{j} Ex_j(t) + \eta_i^{j} Ex_j(t-1)\right) \right)^2 \geq 0, \\
Q_{ii}(t) &=& E\left(\sum_{j=1}^2\left(\sigma_i^{j}\tilde{x}_j(t) + \eta_i^{j}\tilde{x}_j(t-1)\right) \right)^2 \geq 0,
\end{eqnarray*}
and
\begin{eqnarray*}
F_{ij}(t) &=& \sum_{k=1}^2 \int_0^t X_i^k(t-s) X_j^k(t-s) P_{kk}(s) d s, \\
M_{ij}(t) &=& \sum_{k=1}^2 \int_0^t  X_i^k(t-s) X_j^k(t-s) (P_{kk}(s) + Q_{kk}(s)) d s, \\
N_{ij}(t) &=& \sum_{k=1}^2 \int_0^{t-1}  X_i^k(t-s) X_j^k(t-1-s) (P_{kk}(s) + Q_{kk}(s)) d s.
\end{eqnarray*}

Before we state and prove the main results, we introduce some preliminaries below.

When $N=2$, we consider the delay differential equation
\begin{equation}\label{dde2}
\left\{
\begin{array}{rcl}
\displaystyle \dfrac{d x_1}{dt}&=&\displaystyle\sum_{j=1}^2 \left(a_1^j x_j(t) + b_1^j x_j(t-1)\right), \\
\displaystyle \dfrac{d x_2}{dt}&=&\displaystyle\sum_{j=1}^2 \left(a_2^j x_j(t) + b_2^j x_j(t-1)\right).
\end{array}
\right.
\end{equation}
The characteristic function of \eqref{dde2} is given by
\begin{eqnarray}\label{hlambda}
h(\lambda)&=& \det(\Delta(\lambda))=
\det\left( \begin{array}{cc}
\lambda - a_{1}^{1} - b_{1}^{1} e^{-\lambda} & - a_{1}^{2} - b_{1}^{2} e^{-\lambda} \nonumber \\
- a_{2}^{1} - b_{2}^{1} e^{-\lambda} & \lambda - a_{2}^{2} - b_{2}^{2} e^{-\lambda}
\end{array}\right)\\
&=&\lambda^{2}+a\lambda+b+(c \lambda + d)e^{-\lambda} + r e^{-2\lambda},
\end{eqnarray}
where
$$
a = - a_{1}^{1}-a_{2}^{2},\quad b=a_{1}^{1}a_{2}^{2}-a_{1}^{2}a_{2}^{1},\quad c = - b_{1}^{1} - b_{2}^{2},
$$
$$
d = a_{1}^{1} b_{2}^{2} + a_{2}^{2} b_{1}^{1}  - a_{1}^{2} b_{2}^{1} - a_{2}^{1} b_{1}^{2} ,
\quad r = b_{1}^{1} b_{2}^{2}-b_{1}^{2} b_{2}^{1}.
$$
The Laplace transform of the fundamental matrix is
\begin{equation}\label{LX2}
\mathcal{L}(X)(\lambda)=\Delta^{-1}(\lambda)
=\dfrac{1}{h(\lambda)}\left( \begin{array}{cc}
\lambda - a_{2}^{2} - b_{2}^{2} e^{-\lambda} & a_{1}^{2} + b_{1}^{2} e^{-\lambda} \\
a_{2}^{1} + b_{2}^{1} e^{-\lambda} & \lambda - a_{1}^{1} - b_{1}^{1} e^{-\lambda}
\end{array}\right).
\end{equation}

Here we give some properties of the fundamental solution $X(t)$ that are useful for the estimate of the second moment below.

Recall
\begin{equation}
\label{alpha0-2}
\alpha_0 = \sup\{\Re(\lambda): h(\lambda) = 0, \lambda\in \mathbb{C}\}.
\end{equation}
When $\alpha_{0}<0$, from \eqref{LX2}, we have
\begin{equation*}
X(t) = \dfrac{1}{2 \pi i} \lim_{T\to+\infty}\int_{-iT}^{iT}e^{\lambda t}\Delta^{-1}(\lambda) d \lambda = \dfrac{1}{2\pi}\int_{-\infty}^{+\infty}e^{i\omega t}\Delta^{-1}(i\omega) d\omega.
\end{equation*}
Thus, from \eqref{LX2}, we obtain
\begin{eqnarray*}
X_{1}^{1}(t)&=&\dfrac{1}{2\pi}\int_{-\infty}^{+\infty}
\dfrac{i\omega-a_{2}^{2}-b_{2}^{2}e^{-i\omega}}{h(i\omega)}e^{i\omega t} d\omega, \\
X_{1}^{2}(t)&=&\dfrac{1}{2\pi}\int_{-\infty}^{+\infty}
\dfrac{a_{1}^{2}+b_{1}^{2}e^{-i\omega}}{h(i\omega)}e^{i\omega t} d\omega, \\
X_{2}^{1}(t)&=&\dfrac{1}{2\pi}\int_{-\infty}^{+\infty}
\dfrac{a_{2}^{1}+b_{2}^{1}e^{-i\omega}}{h(i\omega)}e^{i\omega t} d\omega, \\
X_{2}^{2}(t)&=&\dfrac{1}{2\pi}\int_{-\infty}^{+\infty}
\dfrac{i\omega-a_{1}^{1}-b_{1}^{1}e^{-i\omega}}{h(i\omega)}e^{i\omega t} d\omega.
\end{eqnarray*}
Obviously, $X_{i}^{k}(t)X_{j}^{k}(t)$ and  $X_{i}^{k}(t)X_{j}^{k}(t-1)\;(i,j,k=1,2)$  have Laplace transforms. When $\alpha_0<0$, explicit
expressions and estimates for the Laplace transforms  
$\mathcal{L}\left(X_{i}^{k}(t)X_{j}^{k}(t)\right)$, $\mathcal{L}\left(X_{i}^{k}(t)X_{j}^{k}(t-1)\right)$
are given in Appendix B.

\subsection{Unboundedness of the second moment}
Here we give a result for the unboundedness of the second moment of \eqref{sdde3} when $\alpha_0 > 0$. First, we note a rare situation
when the coefficients $\mu_i, \sigma_{i}^{j}, \eta_{i}^{j} \,(i,j=1,2)$ satisfy the following assumption.

\begin{minipage}{12cm}
\textbf{Assumption H:}\
\textit{$\mu_{1} = \mu_{2} = 0$, and there is a root $\lambda\in\R$ of $h(\lambda)= 0$ and an eigenvector
$c=\left(\begin{array}{cc}c_{1}\\c_{2}\end{array}\right)\in\R^{2}$ corresponding to the eigenvalue $\lambda$ such that
$$\sum_{j=1}^2(\sigma_i^{j} + e^{-\lambda} \eta_i^{j}) c_{j} =0,\quad i=1,2.
$$}
\end{minipage}

When the assumption H is satisfied, the stochastic equation \eqref{sdde3} has a deterministic solution $x(t)= e^{\lambda t}c\,(t\geq 0)$ with
initial function
$\phi=e^{\lambda t}c\ (-1\leq t<0)$, and the corresponding second moment $M(t) = 0$. This is a very rare situation and is excluded in discussions below.

The following result shows that in general the second moment of \eqref{sdde3} is unbounded when the trivial solution of \eqref{dde} is unstable.

\begin{theorem}\label{Thm-unb}
Let $\alpha_{0}$ be defined as \eqref{alpha0-2} and $\alpha_{0}>0$.
If the assumption H is not satisfied, the second moment of \eqref{sdde3} is unbounded.
\end{theorem}
\Proof
We only need to show that there is a special solution $x(t; \phi)$ of \eqref{sdde3} such
that the corresponding second moment is unbounded. Note that $P_{ii},\ Q_{ii}\geq 0\ (i=1,2)$, and hence
\begin{eqnarray}\label{MXP}
\|M(t)\| & \geq & M_{11}(t) \geq F_{11}(t) \nonumber\\
&=&\sum_{i=1}^2 \int^{t}_{0}X_{1}^{i}(t-s)^2 P_{ii}(s) d s\geq 0.
\end{eqnarray}
Let $\lambda$ be a solution of $h(\lambda)= 0$ with $0 < \Re(\lambda)\leq \alpha_{0}$, and $c\in \mathbb{R}^2$ the corresponding eigenvector.
Then $x_{\phi_1}(t)=\Re(e^{\lambda t}c)$ is the solution of \eqref{dde2} with initial function
$ \phi_1=\Re(e^{\lambda t}c)\ (-1\leq t<0)$.  Since the assumption H is not satisfied, $x_{\phi_1}(t)$ is not a solution of \eqref{sdde3}. Following
the proof of Theorem \ref{Thm-X} (3), there is a subset $U\subset\R^{+}$
with $m(U)=+\infty$ and positive constants $C_{1},\,C_{2}$ such that for any $ t \in U$,
\begin{eqnarray*}
|X_{1}^{i}(t)| \geq  e^{\bar{\alpha}t}, \quad \bar{\alpha} \in(0,\alpha_{0}),\; \;i=1,2
\end{eqnarray*}
and
$$
P_{11}(t)  \geq  C_{1} \quad \text{or} \quad P_{22}(t) \geq C_{2}.
$$
Thus, from \eqref{MXP}, the second moment $M(t)$ is unbounded.
\qed

\subsection{Boundedness of the second moments}

Now, we always assume $\alpha_{0}< 0$. Then \eqref{LM}-\eqref{LQ} become
\begin{eqnarray}
\label{LMij}
\mathcal{L}(M_{ij}) &=& \sum_{k=1}^2\mathcal{L}(X_i^k X_j^k) \left(\mathcal{L}(P_{kk}) + \mathcal{L}(Q_{kk})\right),\\
\label{LNij}
\mathcal{L}(N_{ij}) &=& \sum_{k=1}^2\mathcal{L}(X_i^k(t)X_j^k(t-1)) \left(\mathcal{L}(P_{kk}) + \mathcal{L}(Q_{kk})\right),
\end{eqnarray}
and
\begin{eqnarray}
\label{LQii}
\mathcal{L}(Q_{ii}) &=& \sum_{k=1}^2\left((\sigma_i^{k})^2 + (\eta_{i}^{k})^2 e^{-\lambda} \right)\mathcal{L}(M_{kk})
+ 2 \sum_{k,l=1}^{2}\sigma_{i}^{k} \eta_{i}^{l} \mathcal{L}(N_{kl})\nonumber\\
&&{} + 2\left(\sigma_i^{1} \sigma_i^{2} + \eta_{i}^{1} \eta_{i}^{2}e^{-\lambda}  \right)\mathcal{L}(M_{12}).
\end{eqnarray}

\begin{proposition} \label{Proposition}
Define the matrices
\begin{equation}
\label{AG}
A(\lambda)=\left(\begin{array}{cc}
\mathcal{L}((X_1^1)^2) \;\; & \mathcal{L}((X_1^2)^2)\\[0.1cm]
\mathcal{L}((X_2^1)^2) \;\; & \mathcal{L}((X_2^2)^2)
\end{array}\right),\quad
G(\lambda)=\left(\begin{array}{cc}
G_1^1(\lambda) \;\; & G_1^2(\lambda) \\[0.1cm]
G_2^1(\lambda) \;\; & G_2^2(\lambda)
\end{array}\right)
\end{equation}
with
\begin{eqnarray}
\label{eq:G}
&& G_k^q(\lambda) = (\sigma_k^q)^2 + (\eta_k^q)^2 e^{-\lambda} + \sum_{p=1}^2 T_k^p S_p^q \quad (k,q=1,2),\\
\label{eq:T}
&& T_k^p = 2 (\sigma_k^1 \sigma_k^2 + \eta_k^1 \eta_k^2 e^{-\lambda})\mathcal{L}(X_1^p X_2^p)
+ 2\sum_{m,l=1}^2 \sigma_k^m \eta_k^l \mathcal{L}(X_m^p(t)X_l^p(t-1)) \quad \qquad
\end{eqnarray}
and $S=(S_p^q)_{2\times 2}=A(\lambda)^{-1}$ for $\Re(\lambda)> \alpha_{A}$, where $$
\alpha_{A} \triangleq \sup \left\{\Re(\lambda): \det(A(\lambda)) = 0, \lambda\in \C \right\}.
$$
Then
\begin{equation}\label{LM12}
\mathcal{L}({M}_{12}) = \sum_{k,i=1}^2\mathcal{L}(X_1^kX_2^k)S_k^i\mathcal{L}(M_{ii})
\end{equation}
and
\begin{equation}
\label{LM4}
\left(\begin{array}{c}
\mathcal{L}(M_{11})\\
\mathcal{L}(M_{22})
\end{array}\right)
= (I-D(\lambda))^{-1}A(\lambda)
\left(\begin{array}{c}
\mathcal{L}(P_{11}) \\
\mathcal{L}(P_{22})
\end{array}\right),
\end{equation}
where $D(\lambda) =A(\lambda)G(\lambda)$.
\end{proposition}

\Proof
Assume $\Re(\lambda)>\alpha_{A}$ and therefore $S=A(\lambda)^{-1}$ is well defined. From \eqref{LMij}, we have
$$
\mathcal{L}(M_{ii})= \sum_{k=1}^2\mathcal{L}((X_i^k)^2)
(\mathcal{L}(P_{kk})+\mathcal{L}(Q_{kk})), \quad i=1,2,
$$
therefore
\begin{equation}\label{LPQ}
\mathcal{L}(P_{kk})+\mathcal{L}(Q_{kk}) = \sum_{i=1}^2S_k^i \mathcal{L}(M_{ii}), \quad k=1,2,
\end{equation}
which gives \eqref{LM12} from \eqref{LMij}.

From \eqref{LNij} and \eqref{LPQ},
$$
\mathcal{L}(N_{ij}) = \sum_{p,q=1}^2S_p^q \mathcal{L}(X_i^p(t)X_j^p(t-1))\mathcal{L}(M_{qq}).
$$
Thus, from \eqref{LQii}, \eqref{AG}-\eqref{eq:T} and \eqref{LM12}, we obtain
\begin{eqnarray}
\label{LQii1}
\mathcal{L}(Q_{ii}) &=& \sum_{q=1}^2 \left[((\sigma_i^q)^2 + (\eta_i^q)^2 e^{-\lambda})
+ 2\sum_{m,l,p=1}^2 \sigma_i^m \eta_i^l S_p^q \mathcal{L}(X_m^p(t) X_l^p(t-1))\right] \mathcal{L}(M_{qq})\nonumber\\
&&\quad\quad + 2 (\sigma_i^1 \sigma_i^2 + \eta_i^1 \eta_i^2 e^{-\lambda}) \sum_{p=1}^2\mathcal{L}(X_1^p X_2^p)
\sum_{q=1}^2 S_p^q\mathcal{L}(M_{qq}) \nonumber\\
&=& \sum_{q=1}^2 \bigg[((\sigma_i^q)^2 + (\eta_i^q)^2 e^{-\lambda})
+ \sum_{p=1}^2 S_p^q \Big( 2 (\sigma_i^1 \sigma_i^2 + \eta_i^1 \eta_i^2 e^{-\lambda})\mathcal{L}(X_1^p X_2^p) \nonumber\\
&&\quad \quad \, + 2\sum_{m,l=1}^2 \sigma_i^m \eta_i^l \mathcal{L}(X_m^p(t) X_l^p(t-1)) \Big)\bigg] \mathcal{L}(M_{qq})\nonumber\\
&=& \sum_{q=1}^2  G_i^q  \mathcal{L}(M_{qq}).
\end{eqnarray}
Hence, from \eqref{LMij} and \eqref{LQii1}, we have
\begin{eqnarray*}
\mathcal{L}(M_{ii})&=&\sum_{k=1}^2\mathcal{L}(X_i^k X_i^k)\mathcal{L}(P_{kk})
+ \sum_{k=1}^2 \mathcal{L}(X_i^k X_i^k)\sum_{q=1}^2 G_k^q \mathcal{L}(M_{qq}),
\end{eqnarray*}
i.e.,
\begin{equation*}
(I-D(\lambda))\left(\begin{array}{c}
\mathcal{L}(M_{11})\\
\mathcal{L}(M_{22})
\end{array}\right) = A(\lambda)\left(\begin{array}{c}
\mathcal{L}(P_{11})\\
\mathcal{L}(P_{22})
\end{array}\right),
\end{equation*}
which yields \eqref{LM4}, and the Proposition is proved.
\qed

Denote
\begin{equation}\label{H2}
H(\lambda) = \mathrm{det}(I-D(\lambda)).
\end{equation}
From \eqref{LM4},
\begin{equation}
\label{LMP}
\left(\begin{array}{c}
\mathcal{L}(M_{11})\\
\mathcal{L}(M_{22})
\end{array}\right) = \dfrac{\mathrm{adj}(I-D(\lambda))}{H(\lambda)} A(\lambda)
\left(\begin{array}{c}
\mathcal{L}(P_{11})\\
\mathcal{L}(P_{22})
\end{array}\right),
\end{equation}
where $\mathrm{adj}(\cdot)$ denotes the adjoint matrix.

Let
\begin{equation}
\label{alpha-bar0}
\bar{\alpha}_0 = \sup\{\Re(\lambda): h(\lambda)\det(A(\lambda)) = 0, \lambda\in \mathbb{C}\}.
\end{equation}
Then $\bar{\alpha}_0 = \max\{\alpha_0, \alpha_{A}\}$ \footnote{We conjecture that $\bar{\alpha}_0 = \alpha_0$, but are not able to prove.}, and $A(\lambda)$ is invertible for $\Re(\lambda) > \bar{\alpha}_0$.
Thus $D(\lambda)$ and $\dfrac{\mathrm{adj}(I-D(\lambda))}{H(\lambda)}$ are analytic for $\Re(\lambda) > \bar{\alpha}_0$.

In the following theorem, we show that $H(\lambda)$ is the characteristic function for the second moment boundedness.

\begin{theorem}\label{Thm-bound}
Let $H(\lambda)=\det(I-D(\lambda))$ and $D(\lambda)$ be defined as in Proposition \ref{Proposition}. Let $\bar{\alpha}_{0}$ be defined
at \eqref{alpha-bar0} and assume $\bar{\alpha}_{0}<0$. Then
\begin{enumerate}
\item[\textnormal{(i)}] if all roots of the equation $H(\lambda)=0$ have negative real parts,
then the second moment for any solution of \eqref{sdde3} is bounded, and $M(t)$ approaches a $2\times 2$ constant matrix
exponentially as $t\rightarrow +\infty;$
\item[\textnormal{(ii)}] if the equation $H(\lambda)=0$ has a root with positive real part,
and the assumption H is not satisfied, then there exists a solution of \eqref{sdde3} whose second moment is unbounded.
\end{enumerate}
\end{theorem}

To prove Theorem \ref{Thm-bound}, we first give some useful Lemmas.

\begin{lemma}\label{Lem-D}
Let $D(\lambda)$ be defined as in Proposition \ref{Proposition} and assume $\bar{\alpha}_{0}<0$. Then there exist constants $d_0$ and $T_0$
such that for $|\lambda|\geq T_{0}$ and $\Re(\lambda)>\bar{\alpha}_0$,
\begin{equation}\label{dij2}
\| D(\lambda)\|\leq \dfrac{d_0}{|\lambda|}.
\end{equation}
\end{lemma}
\Proof From Lemma \ref{LemLyij}, there exists $R_0$ so that $\|A(\lambda)\| < \dfrac{R_0}{|\lambda|}$ when $|\lambda|$ is large enough. Thus,
it is enough to show that there exists a constant $g_0$ such that when $\Re(\lambda)>\bar{\alpha}_0$,
\begin{equation}\label{G2}
\limsup_{|\lambda|\to +\infty}\left|G_k^q(\lambda)\right| < g_0,\quad \forall k,q=1,2,
\end{equation}
Then \eqref{dij2} is satisfied with $d_0 = R_0 g_0$.

To prove \eqref{G2}, we only need to show that when $\Re(\lambda)>\bar{\alpha}_0$,
\begin{equation}\label{Spq}
\begin{array}{l}
\displaystyle \limsup_{|\lambda|\to +\infty}\left|S_p^q \mathcal{L}(X_1^p X_2^p) \right| < +\infty,\\
\displaystyle \limsup_{|\lambda|\to +\infty}\left| S_p^q \mathcal{L}(X_m^p(t)X_l^p(t-1))\right| < +\infty.
\end{array}
\end{equation}
Here we only give the proof of the first result of \eqref{Spq} for $p=q=1$, and the others are similar.

For any $\Re(\lambda)>\bar{\alpha}_0,$  from Lemma \ref{LemLyij}, we obtain
\begin{eqnarray*}
S_1^1(\lambda) = \dfrac{\mathcal{L}((X_2^2)^2)}{\det(A(\lambda))}
=\dfrac{\lambda-\alpha_0-a_{1}^{1}-b_{1}^{1}e^{-(\lambda-\alpha_0)}+h(\lambda-\alpha_0)g_{22}(\lambda)}{\det(A(\lambda))},
\end{eqnarray*}
where
\begin{eqnarray*}
\det(A(\lambda)) &=& \mathcal{L}((X_1^1)^2)(\lambda) \mathcal{L}((X_2^2)^2)(\lambda)
-\mathcal{L}((X_1^2)^2)(\lambda) \mathcal{L}((X_2^1)^2)(\lambda)\\
&=& \left[\lambda-\alpha_0-a_{2}^{2}-b_{2}^{2}e^{-(\lambda-\alpha_0)}+ h(\lambda-\alpha_0) g_{11}(\lambda)\right] \\
&& {} \times  \big[\lambda-\alpha_0-a_{1}^{1}-b_{1}^{1}e^{-(\lambda-\alpha_0)} + h(\lambda-\alpha_0) g_{22}(\lambda)\big]  \\
&& {} -\left[a_{1}^{2}+b_{1}^{2}e^{-(\lambda-\alpha_0)}+h(\lambda-\alpha_0)g_{12}(\lambda)\right]
\left[a_{2}^{1}+b_{2}^{1}e^{-(\lambda-\alpha_0)}+h(\lambda-\alpha_0)g_{21}(\lambda)\right].
\end{eqnarray*}
Thus, from \eqref{Xij2}, when $\Re(\lambda)>\bar{\alpha}_0$, $\det(A(\lambda))\not= 0$ and there exists a constant $s_{11}>0$ such that
\begin{equation*}
\lim_{|\lambda|\to +\infty}\left|\lambda S_1^1(\lambda) \right|
=\lim_{|\lambda|\to +\infty} \dfrac{ \left|\lambda \left(\lambda-\alpha_0-a_{1}^{1}-b_{1}^{1}e^{-(\lambda-\alpha_0)}
+h(\lambda-\alpha_0)g_{22}(\lambda)\right)\right|}{|\det(A(\lambda))|}  = s_{11}
\end{equation*}
Hence, from \eqref{Xij3}, there exist two positive constants $\tilde{T}$ and $\tilde{R}_{11}=2s_{11}$ such that for $|\lambda|>\tilde{T}$
and $\Re(\lambda)>\bar{\alpha}_0$,
\begin{equation}\label{eq:S11}
\left|\lambda^2 S_1^1 \mathcal{L}\left(X_1^1 X_2^1\right)\right| \leq \tilde{R}_{11},
\end{equation}
which implies
$$\limsup_{|\lambda|\to\infty}|S_1^1\mathcal{L}(X_1^1X_1^2)| < +\infty,$$
and the first inequality in \eqref{Spq} for $p=q=1$ is proved.
\qed

From Lemma \ref{Lem-D}, when $\Re(\lambda)>\bar{\alpha}_{0}$, $H(\lambda)$ is analytic and
\begin{equation*}
\lim_{|\lambda|\to +\infty}|H(\lambda)|=1.
\end{equation*}
Thus there is a real number $\beta_0$ such that all roots of $H(\lambda)=0$ satisfy
$\Re(\lambda) < \beta_0$ (refer to the discussion in \cite[Lemma 4.1 in Chapter 1]{Hale93}), where
\begin{equation}
\beta_{0}=\sup\{\Re(\lambda): H(\lambda)=0, \lambda\in\C\}.\label{beta0}
\end{equation}

Let
\begin{equation}
Z(t) = \mathcal{L}^{-1}\left(\frac{1}{H(\lambda)}\mathrm{adj}(I-D(\lambda)) A(\lambda)\right),
\end{equation}
then from \eqref{LMP},
\begin{equation}
\label{eq:MP}
\left(
\begin{array}{rcl}
M_{11}(t)\\
M_{22}(t)
\end{array}\right) = \int_0^t Z(t-s) \left(
\begin{array}{rcl}
P_{11}(s)\\
P_{22}(s)
\end{array}\right) ds.
\end{equation}

The following result gives an estimation of $Z(t)$.

\begin{lemma}\label{Lem-Z}
Let $\beta_0$ be defined at \eqref{beta0} and assume $\bar{\alpha}_{0} <0$. There exists a positive constant
$C_3$ such that for any $\beta>\max\{\bar{\alpha}_0, \beta_{0}\}$,
\begin{equation}\label{ZC3}
\|Z(t)\| < C_3 e^{\beta t}.
\end{equation}
\end{lemma}
The proof is similar to that of \cite[Theorem 5.2 in Chapter 1]{Hale93} and the details are given in Appendix C.

\begin{lemma}
\label{LemM12}
Assume $\bar{\alpha}_{0} <0$. If $M_{11}(t), M_{22}(t)$ are bounded, then $M_{12}(t)$ is also bounded for any $t > 0$.
\end{lemma}
\Proof
Assume that there exists a positive constant $M_0$ so that
$$
M_{ii}(t) \leq  M_0,\quad \forall \, t>0, \;i=1,2.
$$
Let
$$
Y_i(t) = \mathcal{L}^{-1}\left(\sum_{k=1}^2\mathcal{L}(X_1^kX_2^k) S_k^i\right),\quad i=1,2.
$$
Then \eqref{LM12} yields
\begin{equation}
\label{M12}
M_{12}(t) = \sum_{i=1}^2 \int_0^t Y_i(t-s) M_{ii}(s) d s.
\end{equation}
Similar to the proof of \eqref{eq:S11}, there exist positive constants $\tilde{R}_{ik}\,(i,k=1,2)$ and $T^*$ such that
for $|\lambda|>T^*$ and $\Re(\lambda)>\bar{\alpha}_0$,
\begin{equation*}
\left|\lambda^2 \mathcal{L}(X_1^kX_2^k) S_k^i\right| \leq \tilde{R}_{ik}.
\end{equation*}
Thus, when $\bar{\alpha}_{0} <0$, similar to the proof of Lemma \ref{Lem-Z}, there exist positive constants $\rho_i\,(i=1,2)$ such that
for any $\alpha >\bar{\alpha}_0$,
\begin{equation}\label{Yi-rho}
\|Y_i(t)\| \leq \rho_i e^{\alpha t},\quad i=1,2.
\end{equation}
Hence, from \eqref{Yi-rho}, when $\bar{\alpha}_{0} <0$, for any $\alpha \in \left(\bar{\alpha}_0, 0\right)$,
\begin{eqnarray*}
|M_{12}(t)| &\leq &  M_0 \int_0^t (|Y_1(s)| + |Y_2(s)|) d s \leq  (\rho_1 + \rho_2) M_0 \int_0^t e^{\alpha s}d s \\
& \leq & \dfrac{2(\rho_1 + \rho_2) M_0}{|\alpha|},
\end{eqnarray*}
i.e., $M_{12}(t)$ is bounded.
\qed

Now we are ready to prove Theorem \ref{Thm-bound}.

\noindent\textit{Proof of Theorem \ref{Thm-bound}.}\ \
Let $H(\lambda)=\det(I-D(\lambda)) $ and $D(\lambda)$ be defined as in Proposition \ref{Proposition}, and $\beta_{0}$ be defined as \eqref{beta0}.
We also assume $\bar{\alpha}_{0}<0$.

From \eqref{Ex}, for any $\alpha \in (\bar{\alpha}_0, 0)$, there exists a constant $K_{2}=K_{2}(\phi)$ $(\phi \in C([-1,0], \R^N) )$
such that
\begin{equation}
P_{ii}(t) \leq (\mu_i)^2 + K_2 e^{\alpha t},\quad\forall \, t\geq 0, \, i=1,2.
\end{equation}
The prove is straightforward from
\begin{eqnarray*}
P_{ii}(t)&=&\left(\mu_i + \sigma_i^{1}Ex_{1}(t)+\eta_i^{1}Ex_{1}(t-1)
+\sigma_i^{2} Ex_{2}(t)+\eta_i^{2}Ex_{2}(t-1)\right)^{2} \\
&\leq & \Big( |\mu_i|+\left(|\sigma_i^{1}|+|\sigma_i^{2}|\right)\|Ex(t)\|
+\left(|\eta_i^{1}|+|\eta_i^{2}|\right)\|Ex(t-1)\| \Big)^{2}\\
&\leq &(\mu_i)^2+e^{\alpha t}\Big[\tilde{K}^2\|\phi\|^2\left(|\sigma_i^{1}|+|\sigma_i^{2}|\right)^2
+\tilde{K}^2\|\phi\|^2e^{-2\alpha}\left(|\eta_i^{1}|+|\eta_i^{2}|\right)^2 \\
&&{}+2\tilde{K}\|\phi\||\mu_i|\left(|\sigma_i^{1}|+|\sigma_i^{2}|\right)
+2\tilde{K}\|\phi\|e^{-\alpha}|\mu_i|\left(|\eta_i^{1}|+|\eta_i^{2}|\right)\\
&&{} +2\tilde{K}^2\|\phi\|^2e^{-\alpha}\left(|\sigma_i^{1}|+|\sigma_i^{2}|\right)
\left(|\eta_i^{1}|+|\eta_i^{2}|\right)\Big].
\end{eqnarray*}

(i). Assume $\beta_0 < 0$. Let
$$
\hat{M}(t) = \left(\begin{array}{c}
M_{11}(t)\\
M_{22}(t)
\end{array}\right),\quad \mu_0 = (\mu_1)^2 + (\mu_2)^2.
$$
From Lemma \ref{Lem-Z}, for any $\beta\in \left(\max\{\bar{\alpha}_0, \beta_{0}\},0\right),\;\alpha\in(\bar{\alpha}_{0},0)$ and $\alpha\not=\beta$,
from \eqref{ZC3},
\begin{eqnarray*}
\|\hat{M}\| &\leq & \int^{t}_{0}\|Z(s)\|(P_{11}(t-s)+ P_{22}(t-s))d s  \\
&\leq &C_{3}\int^{t}_{0}e^{\beta s}\left(\mu_0 +2 K_2e^{\alpha (t-s)}\right)d s \\
&=&C_{3}\left( \mu_0 \dfrac{1-e^{\beta t}}{-\beta} +2 K_2 \dfrac{e^{\alpha t} - e^{\beta t}}{\alpha-\beta}\right) \\
&\leq &C_3 \left(\dfrac{\mu_0}{|\beta|} + \dfrac{4 K_2}{|\alpha- \beta|}\right).
\end{eqnarray*}
From Lemma \ref{LemM12}, the second moment $M(t)$ is bounded for any initial function $\phi\in C([-1,0], \R^2).$

Let
\begin{equation*}
\hat{M}_{\infty}=\left(\begin{array}{c}
\hat{M}_{\infty,1} \\
\hat{M}_{\infty,2}
\end{array}\right)
=
\int^{+\infty}_{0} Z(s) \left(\begin{array}{c}
(\mu_1)^2\\
(\mu_2)^2
\end{array}\right)ds.
\end{equation*}
It is easy to have
\begin{equation}\label{hat-M}
\|\hat{M}_{\infty}\| \leq  \mu_0 \int^{+\infty}_{0} \|Z(s)\| d s
\leq  C_{3}\mu_0 \int^{+\infty}_{0}e^{\beta s}ds
= \dfrac{C_{3}\mu_0 }{|\beta|}.
\end{equation}
Thus, from \eqref{ZC3},
\begin{eqnarray}
\left\|\hat{M}(t)-\hat{M}_{\infty}\right\|
&=& \left\|\int^{t}_{0} Z(s)\left(
\begin{array}{c}
P_{11}(t-s)-(\mu_1)^2\\
P_{22}(t-s)-(\mu_2)^2
\end{array}\right) d s
-\int^{+\infty}_{t}Z(s) \left(\begin{array}{c}
(\mu_1)^2\\
(\mu_2)^2
\end{array}\right)ds\right\|\nonumber\\
&\leq & 2 K_2\int^{t}_{0} \|Z(s)\| e^{\alpha(t-s)}d s  + \mu_0 \int^{+\infty}_{t} \|Z(s)\| d s\nonumber\\
&\leq & C_{3}\mu_0 \dfrac{e^{\beta t}}{-\beta}
+2 C_{3} K_2 \dfrac{e^{\alpha t}-e^{\beta t}}{\alpha-\beta}\nonumber\\
&\leq & C_3\left(\dfrac{\mu_0}{|\beta|}
+ \dfrac{4 K_2}{|\alpha-\beta|}\right)e^{t\max\{\alpha, \beta\}}
\triangleq C_{4}e^{t\max\{\alpha, \beta\}},
\label{eq:C4}
\end{eqnarray}
which implies $\left\|\hat{M}(t)-\hat{M}_{\infty}\right\|\to 0\,( \text{as}\,t\to+\infty)$ exponentially  since
$\max\{\alpha, \beta\}<0$.

From \eqref{M12}, and let
$$
M_\infty = \left(
\begin{array}{cc}
\hat{M}_{\infty,1} & \displaystyle \sum_{i=1}^2 \int_0^{+\infty} Y_i(s) \hat{M}_{\infty,i} d s\\
\displaystyle \sum_{i=1}^2 \int_0^{+\infty} Y_i(s) \hat{M}_{\infty,i} d s  & \hat{M}_{\infty,2}
\end{array}
\right).
$$
Obviously, from \eqref{Yi-rho} and \eqref{hat-M}, $M_\infty$ is bounded.
Thus
\begin{eqnarray*}
\left\|M(t)-M_{\infty}\right\|
&=&  \left(
\begin{array}{cc}
M_{11}-\hat{M}_{\infty,1} & \displaystyle M_{12}-\sum_{i=1}^2 \int_0^{+\infty} Y_i(s) \hat{M}_{\infty,i} d s\\
\displaystyle M_{12}-\sum_{i=1}^2 \int_0^{+\infty} Y_i(s) \hat{M}_{\infty,i} d s  & M_{22}-\hat{M}_{\infty,2}
\end{array}\right)\\
&=& \left\|\hat{M}(t)-\hat{M}_{\infty}\right\| + 2 \left| M_{12}-\sum_{i=1}^2 \int_0^{+\infty} Y_i(s) \hat{M}_{\infty,i} d s \right|.
\end{eqnarray*}
From \eqref{M12}, \eqref{Yi-rho}, \eqref{hat-M} and \eqref{eq:C4}, we obtain for any $\tilde{\alpha}\in(\bar{\alpha}_{0},0)$ and
$\tilde{\alpha} \not= \max\{\alpha, \beta\}$,
\begin{eqnarray*}
&& \left| M_{12}-\sum_{i=1}^2 \int_0^{+\infty} Y_i(s) \hat{M}_{\infty,i} d s \right| \\
& \leq & \sum_{i=1}^2 \left(\int_0^{t} \left|Y_i(s)\right| \left|M_{ii}(t-s)- \hat{M}_{\infty,i}\right| d s
+ \int_t^{+\infty} \left|Y_i(s)\right| \left|\hat{M}_{\infty,i}\right| d s \right) \\
& \leq &  C_4 (\rho_1 + \rho_2) \int_0^{t}  e^{\tilde{\alpha} s} e^{(t-s)\max\{\alpha, \beta \}} d s
+ \dfrac{C_{3}\mu_0 (\rho_1 + \rho_2)}{|\beta|} \int_t^{+\infty}  e^{\tilde{\alpha} s} d s  \\
& = &  C_4 (\rho_1 + \rho_2) \dfrac{ e^{\tilde{\alpha} t} - e^{t\max\{\alpha, \beta \}}}{\tilde{\alpha}-\max\{\alpha, \beta \}}
+ \dfrac{C_{3}\mu_0 (\rho_1 + \rho_2)}{|\tilde{\alpha}\beta|}e^{\tilde{\alpha} t}  \\
& \leq & \left( \dfrac{ 2C_4 (\rho_1 + \rho_2) }{|\tilde{\alpha}-\max\{\alpha, \beta \}|}
+ \dfrac{C_{3}\mu_0 (\rho_1 + \rho_2)}{|\tilde{\alpha}\beta|} \right)e^{t\max\{\tilde{\alpha}, \alpha, \beta \}}.
\end{eqnarray*}
Therefore
\begin{eqnarray*}
\left\|M(t)-M_{\infty}\right\|
&\leq &   \left( C_{4} + \dfrac{ 4C_4 (\rho_1 + \rho_2) }{|\tilde{\alpha}-\max\{\alpha, \beta \}|}
+ \dfrac{ 2C_{3}\mu_0 (\rho_1 + \rho_2)}{|\tilde{\alpha}\beta|} \right)e^{t\max\{\tilde{\alpha}, \alpha, \beta\}},
\end{eqnarray*}
which implies that $M(t)$ approaches to $M_{\infty}$ exponentially as $t\to+\infty$.

(ii). Now we assume $\beta_{0}>0.$ We only need to show that there is a special solution $x(t; \phi)$ such that the corresponding second moment
is unbounded. Similar to the proof of Theorem \ref{Thm-unb}, let $\lambda=\alpha+i \omega\,(\alpha\leq \alpha_{0}<0)$ be a solution of
$h(\lambda)= 0$, and $c=\left(\begin{array}{cc}c_{1}\\c_{2}\end{array}\right)\in\R^{2}$ is an eigenvector corresponding to the eigenvalue $\lambda$,
then $x_{\phi_2}(t)=\Re(e^{\lambda t}c)\,(t\geq -1)$ is a solution of \eqref{dde} with initial function $\phi_2=\Re(e^{\lambda t}c)\in C([-1,0],\R^{2})$.
Hence, for this particular initial function $\phi_2$, since the assumption H is not satisfied, $x_{\phi_2}(t)$ is not a solution of \eqref{sdde2}
and therefore $P_{11}(t)$ or $P_{22}(t)$ is nonzero. Thus the Laplacian $\mathcal{L}(P_{11})$ or $\mathcal{L}(P_{22})$ is nonzero.

Since $\|M(t)\| \geq \|\hat{M}(t)\|\geq M_{11}(t)$, in the following, we only need to show that $M_{11}(t)$
is unbounded for the initial function $\phi_2$. From \eqref{LMP}, we have
\begin{eqnarray*}
M_{11}(t)&=&\dfrac{1}{2\pi i}\lim_{T\to+\infty}
\int_{\bar{c} - i T}^{\bar{c} + i T}e^{st}\Big[\dfrac{(1-d_{22})\mathcal{L}((X_1^1)^2)}{H(s)} \mathcal{L}(P_{11})(s)\\
&&\qquad \qquad \qquad \qquad \quad \;
+ \dfrac{d_{12}\mathcal{L}((X_2^1)^2)}{H(s)} \mathcal{L}(P_{22})(s)\Big]d s,
\end{eqnarray*}
where $\bar{c}>\beta_{0}$.
Here $\dfrac{(1-d_{22})\mathcal{L}((X_1^1)^2)}{H(s)}$ and $\dfrac{d_{12}\mathcal{L}((X_2^1)^2)}{H(s)}$ are analytic functions for
$\Re(s)>\bar{\alpha}_0$ and are nonzeros. It is easy to see that $\mathcal{L}(P_{11})(s),\,\mathcal{L}(P_{22})(s)$ are analytic for
$\Re(s) = \bar{c} > 0$. Thus, similar to the proof of Theorem \ref{Thm-X} (3), there is a sequence $\{t_{k}\}_{k\geq 1}$ such that
$t_{k}\to+\infty$ and $M_{11}(t_{k})\to+\infty$ as $k\to+\infty$, which implies that the second moment is unbounded. Thus, the theorem is proved.
\qquad$\square$

The critical case of $\bar{\alpha}_{0}<0,\ \beta_{0}=0$ is not considered here,  and the issue of boundedness criteria remains open.

The characteristic  function $H(\lambda)$ depends  not only  on the coefficients
of equation \eqref{sdde3}, but also on the Laplace transforms of  $X_{i}^{k}(t)X_{j}^{k}(t)$
and $X_{i}^{k}(t)X_{j}^{k}(t-1)\;(i,j,k=1,2)$. One can calculate these functions numerically according to
Lemma \ref{LemLyij} and \ref{LemLyij2}, however, it is difficult to obtain $\beta_0 = \sup\{\Re(\lambda): H(\lambda) = 0\}$
for a given equation. Hence the sufficient conditions for the second moment boundedness of \eqref{sdde3} established in Theorem \ref{Thm-bound}
are not practical. In applications, one need to derive useful criteria in terms of equation coefficients, either from the proposed characteristic function or following the procedure in the above discussions. Here,
for applications, we give a practical condition for the boundedness of the second moment.

\begin{theorem}\label{Thm-BK}
Assume $\alpha_0 <0$. If there exists $\alpha\in(\alpha_{0},0)$ and the a positive constant $\bar{K}=\bar{K}(\alpha)$ so that
$$
 \|X(t)\| \leq \bar{K} e^{-\alpha t}, \;t>0
$$ and
\begin{equation}
\big(|\sigma_1^{1}|+|\eta_1^{1}|+|\sigma_1^{2}|+|\eta_1^{2}|\big)^{2}+
\big(|\sigma_2^{1}|+|\eta_2^{1}|+|\sigma_2^{2}|+|\eta_2^{2}|\big)^{2}
< -\dfrac{\alpha}{2\bar{K}^{2}},\label{sig-eta}
\end{equation}
then the second moment $M(t)$ is bounded.
\end{theorem}
\Proof
From the expression of $Q_{ii}(t)$ and \eqref{Nij}, we obtain for $i=1,2$,
\begin{eqnarray*}
Q_{ii}(t)&\leq&
\big(|\sigma_i^{1}|+|\sigma_i^{2}|\big)\big(|\sigma_i^{1}|+|\eta_i^{1}|+|\sigma_i^{2}|+|\eta_i^{2}|\big)
\big[M_{11}(t)+2|M_{12}(t)|+M_{22}(t)\big]\\
&& +\big(|\eta_i^{1}|+|\eta_i^{2}|\big)\big(|\sigma_i^{1}|+|\eta_i^{1}|+|\sigma_i^{2}|+|\eta_i^{2}|\big)\big[ M_{11}(t-1)\\
&&\quad +2|M_{12}(t-1)|+M_{22}(t-1)\big]\\
&=&\big(|\sigma_i^{1}|+|\sigma_i^{2}|\big)\big(|\sigma_i^{1}|+|\eta_i^{1}| +|\sigma_i^{2}|+|\eta_i^{2}|\big)\|M(t)\| \\
&& +\big(|\eta_i^{1}|+|\eta_i^{2}|\big)\big(|\sigma_i^{1}|+|\eta_i^{1}|
+|\sigma_i^{2}|+|\eta_i^{2}|\big)\|M(t-1)\|.
\end{eqnarray*}
Since $Q_{ii}(t)\geq 0$, for any $\alpha\in(\alpha_{0},0)$ and $i,j=1,2$,
\begin{eqnarray*}
|M_{ij}(t)| &\leq &|F_{ij}(t)| + \sum_{k=1}^2\int^{t}_{0}|X_i^k(t-s)X_j^k(t-s)Q_{kk}(s)|ds\\
&\leq & |F_{ij}(t)| + \alpha_{2} \bar{K}^{2} e^{2\alpha t}\int^{t}_{0}e^{-2\alpha s}\|M(s)\|ds \\
&&\quad + \alpha_{3} \bar{K}^{2}e^{2\alpha t}\int^{t}_{0}e^{-2\alpha s}\|M(s-1)\|ds,
\end{eqnarray*}
where
\begin{eqnarray*}
\alpha_{2}&=&\big(|\sigma_1^{1}|+|\sigma_1^{2}|\big)\big(|\sigma_1^{1}|+|\eta_1^{1}| +|\sigma_1^{2}|+|\eta_1^{2}|\big)\\
&&{} +\big(|\sigma_2^{1}|+|\sigma_2^{2}|\big)\big(|\sigma_2^{1}|+|\eta_2^{1}| +|\sigma_2^{2}|+|\eta_2^{2}|\big),\\[0.1cm]
\alpha_{3}&=&\big(|\eta_1^{1}|+|\eta_1^{2}|\big)\big(|\sigma_1^{1}|+|\eta_1^{1}| +|\sigma_1^{2}|+|\eta_1^{2}|\big)\\
&&{} +\big(|\eta_2^{1}|+|\eta_2^{2}|\big)\big(|\sigma_2^{1}|+|\eta_2^{1}| +|\sigma_2^{2}|+|\eta_2^{2}|\big).
\end{eqnarray*}
Thus from Lemma \ref{Lem-FL}, we have
\begin{eqnarray*}
\|M(t)\| =\sum^{2}_{i,j=1}|M_{ij}|
&\leq &  K_{1}(1-e^{2\alpha t}) +4 \alpha_{2} \bar{K}^{2} e^{2\alpha t}
\int^{t}_{0}e^{-2\alpha s}\|M(s)\| d s\\
&&{} + 4 \alpha_{3} \bar{K}^{2}e^{2\alpha t} \int^{t}_{0}e^{-2\alpha s}\|M(s-1)\| d s.
\end{eqnarray*}

Let
\begin{equation*}
y(t)=e^{-2\alpha t}\|M(t)\|,\quad  r(t)= K_{1}(e^{-2\alpha t}-1)
\end{equation*}
and
\begin{equation*}
\bar{p} =4 \alpha_{2} \bar{K}^{2}, \quad \bar{q} =4 \alpha_{3}\bar{K}^{2}e^{-2\alpha}.
\end{equation*}
Then
\begin{equation*}
y(t)\leq \bar{p} \int_{0}^{t}{y(s)ds}+ \bar{q} \int_{0}^{t}{y(s-1)ds}+r(t),\; t\geq 0.
\end{equation*}
Choose $\lambda=-2\alpha>0,$ we get
\begin{equation*}
\sup_{t\geq0}|r(t)e^{-\lambda t}| =\sup_{t\geq0} K_{1} (1-e^{2\alpha t}) \leq 2 K_{1} <+\infty.
\end{equation*}
By \eqref{sig-eta}, we obtain
\begin{equation*}
\lambda- \bar{p} - \bar{q} e^{-\lambda} =-2\alpha-4(\alpha_{2}+\alpha_{3})\bar{K}^{2} >-2\alpha+4\bar{K}^{2}
\dfrac{\alpha}{2\bar{K}^{2}}=0.
\end{equation*}
Therefore from Lemma 3.9 in \cite{LeiM}, there exists $C_{5}=C_{5}(\alpha)$ such that
\begin{equation*}
\|M(t)\| e^{-2\alpha t}=y(t)\leq C_{5}e^{-2\alpha t}, \quad t\geq 0,
\end{equation*}
that is, $\|M(t)\| \leq C_{5}$ for all $t\geq 0.$
\qed

In this paper, we have established framework procedure to calculate the characteristic function for the second moment boundedness of linear delay differential equations with a single discrete delay, we also applied the procedure to study a special case of $2$-dimensional equations. However, as we have seen, the resulting function has a very complicate form. These complicate results is in fact show the elaborate correlations of non-Markov processes when both delay and stochastic effects are taken into involved. In spite of the complicate form of final formulations, the procedure of calculating is simple and easy to follow. Thus, in applications, one can develop the characteristic function, for particular equation of studied, following the scheme given here.  We leave these further applications to future works.



\begin{appendix}
\section{Proof of Theorem \ref{Thm-X} (iii)}
\Proof
Let $\alpha_{1}<\alpha_{0}$. Since the zeros of $h(\lambda)$ are isolated,
we can take $\bar{\alpha}\in (\alpha_{1},\alpha_{0})$ such that $\Re(\lambda)=\bar{\alpha}$ does not
contain any root of $h(\lambda)=0$. Next, choose $c_{1}>\alpha_{0}$ and $T>0$, then
\begin{equation}\label{X}
X(t)=\dfrac{1}{2\pi i}\lim_{T\to+\infty}\int_{c_{1}-iT}^{c_{1}+iT}e^{\lambda t}\Delta^{-1}(\lambda)d\lambda.
\end{equation}

To calculate the integral \eqref{X}, we consider the integration of the matrix
$e^{\lambda t}\Delta^{-1}(\lambda)$ around the bounder of the box in the complex plane with
boundary $\gamma = \gamma_{1}\gamma_{2}\gamma_{3}\gamma_{4}$ in the anticlockwise direction, where the segment
$\gamma_{1}$ is the set $\{c_{1}+i\tau: -T\leq \tau \leq T\},$ the segment $\gamma_{3}$ is the set
$\{\bar{\alpha}+i\tau: -T\leq \tau \leq T\}$, the segment $\gamma_{2}$ is the set
$\{u+iT: \bar{\alpha}\leq u \leq c_{1}\}$ and the segment $\gamma_{4}$ is the set $\{u-iT: \bar{\alpha}\leq u \leq c_{1}\}$.
Then Cauchy theorem of residues implies
\begin{equation*}
\oint_{\gamma}e^{\lambda t}\Delta^{-1}(\lambda)d\lambda=2\pi i\sum_{j=1}^{m}
\underset{\lambda=\lambda_{j}}{\Res} e^{\lambda t}\Delta^{-1}(\lambda)\neq 0,
\end{equation*}
where $\lambda_{1}, \lambda_{2}, \cdots, \lambda_{m}$ are roots of $h(\lambda)=0$ inside $\gamma $
($m\geq1$ from the definition of $\alpha_{0}$, and $m<+\infty$ since $h(\lambda)$ is an analytic function).
We also assume that
$$
\bar{\alpha}<\Re(\lambda_{1})\leq \Re(\lambda_{2})\leq \cdots \leq \Re(\lambda_{m})=\alpha_{0}.
$$
Note that
\begin{equation*}
\underset{\lambda=\lambda_{j}}{\Res}e^{\lambda t}\Delta^{-1}(\lambda)=P_{j}(t)e^{\lambda_{j}t},
\end{equation*}
where $P_{j}(t)=\left( P^{kl}_{j}(t) \right)_{N \times N}$ with $P^{kl}_{j}(t)\,(k,l=1, \cdots,N)$ a polynomial of $t$ with degree given by the multiplicity of $\lambda_{j}$ minus $1$.  Thus,
\begin{equation}\label{Pj}
\oint_{\gamma}e^{\lambda t}\Delta^{-1}(\lambda)d\lambda=2\pi i\sum_{j=1}^{m}P_{j}(t)e^{\lambda_{j} t}.
\end{equation}
From the definition of the adjoint matrix, we obtain
\begin{eqnarray*}
\left\|\int_{\gamma_{2}} e^{\lambda t} \Delta^{-1}(\lambda) d\lambda \right\|
& = & \left\|-\int_{\bar{\alpha} + i T}^{c_{1}+ i T} e^{\lambda t} \dfrac{\mathrm{adj}(\Delta(\lambda))}{h(\lambda)} d\lambda\right\| \\
&\leq & e^{c_1 t} \int_{\bar{ \alpha} }^{ c_{1} }  \left\|\dfrac{\mathrm{adj}(\Delta(u+ i T))}{h(u+ i T)}\right\| d u
\to  0 \;( \text{as} \;T\to +\infty),
\end{eqnarray*}
where $\mathrm{adj}(\Delta(\lambda))$ is the adjoint matrix of $\Delta(\lambda)$. In the same way, we have
\begin{equation*}
\left\|\int_{\gamma_{4}}e^{\lambda t}\Delta^{-1}(\lambda)d\lambda\right\|  \to 0 \;( \text{as}\; T \to + \infty).
\end{equation*}
Therefore by \eqref{Pj} we get
\begin{equation*}
\dfrac{1}{2\pi i}\lim_{T\to+\infty}\int_{c_{1}-iT}^{c_{1}+iT}e^{\lambda t}\Delta^{-1}(\lambda)d\lambda
+\dfrac{1}{2\pi i}\lim_{T\to+\infty}
\int_{\bar{\alpha}+iT}^{\bar{\alpha}-iT}e^{\lambda t}\Delta^{-1}(\lambda)d\lambda
=\sum_{j=1}^{m}P_{j}(t)e^{\lambda_{j}t},
\end{equation*}
i.e.,
\begin{equation*}
X(t)=X_{\bar{\alpha}}(t)+\sum_{j=1}^{m}P_{j}(t)e^{\lambda_{j}t},
\end{equation*}
where
\begin{equation*}
X_{\bar{\alpha}}(t)=\dfrac{1}{2\pi i}\lim_{T\to+\infty}
\int_{\bar{\alpha}-iT}^{\bar{\alpha}+iT}e^{\lambda t}\Delta^{-1}(\lambda)d\lambda
\triangleq \left(X^{kl}_{\bar{\alpha}}(t) \right)_{N \times N}.
\end{equation*}
Moreover, similar to the proof of Theorem \ref{Thm-X} (i) (refer \cite{LeiM} or \cite{WLL}), there exists a positive constant
$\bar{C}=\bar{C}(\bar{\alpha})$ such that $X_{\bar{\alpha}}(t)$ satisfies
\begin{equation*}
\left\|X_{\bar{\alpha}}(t)\right\| \leq \bar{C}e^{\bar{\alpha} t}, \quad t\geq 0.
\end{equation*}
Thus we obtain for and $k,l=1, \cdots, N$
\begin{eqnarray*}
|X_{k}^{l}(t)| &\geq &  \left|\sum_{j=1}^{m}P^{kl}_{j}(t)e^{\lambda_{j}t}\right|
-\left|X^{kl}_{\bar{\alpha}}(t)\right|
\geq  \left|\sum_{j=1}^{m}P^{kl}_{j}(t)e^{\lambda_{j}t}\right|-\bar{C}e^{\bar{\alpha}t} \\
&=& e^{\bar{\alpha}t}\left(e^{(\Re(\lambda_{1})-\bar{\alpha})t}f(t)-\bar{C} \right),
\end{eqnarray*}
where
$f(t)= \left|\sum_{j=1}^{m}P^{kl}_{j}(t)e^{(\lambda_{j}-\Re(\lambda_{1}))t}\right|.$

Let $\lambda_{j}=\beta_{j}+i\omega_{j}\,(j=1,2, \cdots, m),$ and assume $n_{0}$ such that $\beta_{j}<\beta_{m}$
when $1\leq j\leq n_{0}$ and $\beta_{j}=\beta_{m}$ when $n_{0}+1\leq j\leq m.$ Then
\begin{eqnarray*}
f(t)&=& e^{(\beta_{m}-\beta_{1})t}
\left|\sum_{j=1}^{n_{0}}e^{-(\beta_{m}-\beta_{j})t} P^{kl}_{j}(t)e^{i\omega_{j}t}
+\sum_{j=n_{0} + 1}^{m} P^{kl}_{j}(t)e^{i\omega_{j}t}\right|\\
&\geq & \left|\sum_{j=n_{0}+1}^{m} \Re\left(P^{kl}_{j}(t)e^{i\omega_{j}t}\right)\right|
-\sum_{j=1}^{n_{0}}e^{-(\beta_{m}-\beta_{j})t} \left|P^{kl}_{j}(t)\right|.
\end{eqnarray*}
Since $P^{kl}_{j}(t)\,(j=1,2, \cdots, m)$ are nonzero polynomials, we find that
\begin{eqnarray*}
\sum_{j=n_{0}+1}^{m} \Re\left(P^{kl}_{j}(t)e^{i\omega_{j}t}\right)
&=&t^{n_{1}}\left[ \sum_{j=n_{0}+1}^{m}a_{j}\cos(\omega_{j}t)+b_{j}\sin(\omega_{j}t)+O(t^{-1}) \right]\\
&=&t^{n_{1}}\left[ \sum_{j=n_{0}+1}^{m}\sqrt{a^{2}_{j}+ b^{2}_{j}}\sin(\theta_{j}+\omega_{j}t)+O(t^{-1}) \right]
\end{eqnarray*}
$(\text{as} \,t\to +\infty),$ where $a_{j},\, b_{j}\,(j=n_{0}+1,\cdots,m)$ are constants and
$a^{2}_{j}+ b^{2}_{j}\neq 0,$ $n_{1}$ is the highest degree of the polynomials $P^{kl}_{j}(t)\,(j=n_{0}+1,\cdots,m)$,
and
$$\sin(\theta_{j})=\dfrac{a_{j}}{\sqrt{a^{2}_{j}+ b^{2}_{j}}},\quad
\cos(\theta_{j})=\dfrac{b_{j}}{\sqrt{a^{2}_{j}+ b^{2}_{j}}}\;\;(j=n_{0}+1,\cdots,m).
$$
Thus, we can always find a
subset $U_{0}\subset \R^+$ with measure $m(U_{0}) = +\infty$ such that all functions
$$
\sin(\theta_{j}+\omega_{j}t)>\varepsilon, \; t\in U_{0}, \; n_{0}+1\leq j \leq m
$$
for some small positive constant $\varepsilon$, and therefore the subset $U$ is
always possible by taking $U= U_{0}\bigcap(t_{0},+\infty)$ with $t_{0}$ large enough. Hence for the
above $\varepsilon$ and $\forall t\in U,$
\begin{equation*}
\left|\sum_{j=n_{0}+1}^{m} \Re\left(P^{kl}_{j}(t)e^{i\omega_{j}t}\right)\right|>2\varepsilon.
\end{equation*}
Furthermore, since
$\sum_{j=1}^{n_0}e^{-(\beta_{m}-\beta_{j})t} \left|P^{kl}_{j}(t)\right|\to 0$ as $t\to+\infty,$
we can take $U$ such that
\begin{equation*}
e^{(\Re(\lambda_{1})-\bar{\alpha})t}f(t)-\bar{C}>1,\quad \forall t\in U
\end{equation*}
and hence for any $ k,l=1,\cdots, N$ and $ t\in U$,
$$
|X_{k}^{l}(t)| \geq e^{\bar{\alpha}t} \left(e^{(\Re(\lambda_{1})-\bar{\alpha})t}f(t)-\bar{C}\right)>  e^{\bar{\alpha}t},
$$
therefore \eqref{Xij} is concluded.
\qed

\section{Expressions and estimates of the Laplace transforms $\mathcal{L}\left(X_{i}^{k}(t) X_{j}^{k}(t)\right)$ and
$\mathcal{L}\left(X_{i}^{k}(t) X_{j}^{k}(t-1)\right)$ }

\begin{lemma}
\label{LemLyij}
Assume $\alpha_0 <0$. For any $\Re(\lambda)>\alpha_0\,(\lambda\in\C),$
\begin{eqnarray}
\label{Xij1-1}
\mathcal{L}((X_1^1)^2)(\lambda)&=&\dfrac{1}{2\pi}\int_{-\infty}^{+\infty}
\dfrac{i\omega-a_{2}^{2}-b_{2}^{2}e^{-i\omega}}{h(i\omega)}
\cdot\dfrac{\lambda-i\omega-a_{2}^{2}-b_{2}^{2}e^{-(\lambda-i\omega)}}{h(\lambda-i\omega)}e^{i\omega t} d\omega \nonumber\\
&=&
\dfrac{\lambda-\alpha_0-a_{2}^{2}-b_{2}^{2}e^{-(\lambda-\alpha_0)}}{h(\lambda-\alpha_0)}+g_{11}(\lambda),\\
\mathcal{L}((X_1^2)^2)(\lambda)&=&\dfrac{1}{2\pi}\int_{-\infty}^{+\infty}
\dfrac{a_{1}^{2}+b_{1}^{2}e^{-i\omega}}{h(i\omega)}
\cdot\dfrac{a_{1}^{2}+b_{1}^{2}e^{-(\lambda-i\omega)}}{h(\lambda-i\omega)} d\omega \nonumber\\
&=&
\dfrac{a_{1}^{2}+b_{1}^{2}e^{-(\lambda-\alpha_0)}}{h(\lambda-\alpha_0)}+g_{12}(\lambda),\\
\mathcal{L}((X_2^1)^2)(\lambda)&=&\dfrac{1}{2\pi}\int_{-\infty}^{+\infty}
\dfrac{a_{2}^{1}+b_{2}^{1}e^{-i\omega}}{h(i\omega)}
\cdot\dfrac{a_{2}^{1}+b_{2}^{1}e^{-(\lambda-i\omega)}}{h(\lambda-i\omega)}d\omega \nonumber\\
&=&
\dfrac{a_{2}^{1}+b_{2}^{1}e^{-(\lambda-\alpha_0)}}{h(\lambda-\alpha_0)}+g_{21}(\lambda),\\
\mathcal{L}((X_2^2)^2)(\lambda)&=&\dfrac{1}{2\pi}\int_{-\infty}^{+\infty}
\dfrac{i\omega-a_{1}^{1}-b_{1}^{1}e^{-i\omega}}{h(i\omega)}
\cdot\dfrac{\lambda-i\omega-a_{1}^{1}-b_{1}^{1}e^{-(\lambda-i\omega)}}{h(\lambda-i\omega)} d\omega \nonumber\\
&=&
\dfrac{\lambda-\alpha_0-a_{1}^{1}-b_{1}^{1}e^{-(\lambda-\alpha_0)}}{h(\lambda-\alpha_0)}+g_{22}(\lambda),\\
\mathcal{L}(X_1^1X_2^1)(\lambda)&=&\dfrac{1}{2\pi}\int_{-\infty}^{+\infty}
\dfrac{a_{2}^{1}+b_{2}^{1}e^{-i\omega}}{h(i\omega)}
\cdot\dfrac{\lambda-i\omega-a_{2}^{2}-b_{2}^{2}e^{-(\lambda-i\omega)}}{h(\lambda-i\omega)}d\omega \nonumber\\
&=&
\dfrac{a_{2}^{1}+b_{2}^{1}e^{-(\lambda-\alpha_0)}}{h(\lambda-\alpha_0)}+g_{31}(\lambda),\\
\label{Xij1-2}
\mathcal{L}(X_1^2X_2^2)(\lambda)&=&\dfrac{1}{2\pi}\int_{-\infty}^{+\infty}
\dfrac{a_{1}^{2}+b_{1}^{2}e^{-i\omega}}{h(i\omega)}
\cdot\dfrac{\lambda-i\omega-a_{1}^{1}-b_{1}^{1}e^{-(\lambda-i\omega)}}{h(\lambda-i\omega)} d\omega \nonumber\\
&=&
\dfrac{a_{1}^{2}+b_{1}^{2}e^{-(\lambda-\alpha_0)}}{h(\lambda-\alpha_0)}+g_{32}(\lambda),
\end{eqnarray}
where
\begin{equation*}
\scriptsize{
\left.\begin{array}{ll}
g_{11}(\lambda)=\dfrac{1}{2\pi}\dint_{-\infty}^{+\infty}
\dfrac{i\omega-a_{2}^{2}-b_{2}^{2}e^{-i\omega}}{h\left(i\omega\right)}\times & \\[0.2cm]
\qquad \qquad
\dfrac{-(a_{2}^{2}+a+\alpha_0)(\lambda-i\omega)-b+\alpha_0 a_{2}^{2}
-\left[(c+b_{2}^{2})(\lambda-i\omega)-\alpha_0 b_{2}^{2}+d\right]e^{-(\lambda-i\omega)}-r e^{-2(\lambda-i\omega)}}
{h\left(\lambda-i\omega\right) \left(\lambda - i\omega -\alpha_0 \right)}d \omega, & \\[0.2cm]
g_{12}(\lambda)=\dfrac{1}{2\pi}\dint_{-\infty}^{+\infty}
\dfrac{a_{1}^{2}+b_{1}^{2}e^{-i\omega}}{h\left(i\omega\right)}\times & \\[0.2cm]
\qquad \qquad
\dfrac{-(\lambda-i\omega)^2+(a_{1}^{2}-a)(\lambda-i\omega)-b-\alpha_0 a_{1}^{2}
-\left[(c-b_{1}^{2})(\lambda-i\omega)+\alpha_0 b_{1}^{2}+d\right]e^{-(\lambda-i\omega)}-r e^{-2(\lambda-i\omega)}}
{h\left(\lambda-i\omega\right) \left(\lambda - i\omega -\alpha_0 \right)}d \omega, & \\[0.2cm]
g_{21}(\lambda)=\dfrac{1}{2\pi}\dint_{-\infty}^{+\infty}
\dfrac{a_{2}^{1}+b_{2}^{1}e^{-i\omega}}{h\left(i\omega\right)}\times & \\[0.2cm]
\qquad \qquad
\dfrac{-(\lambda-i\omega)^2+(a_{2}^{1}-a)(\lambda-i\omega)-b-\alpha_0 a_{2}^{1}
-\left[(c-b_{2}^{1})(\lambda-i\omega)+\alpha_0 b_{2}^{1}+d\right]e^{-(\lambda-i\omega)}-r e^{-2(\lambda-i\omega)}}
{h\left(\lambda-i\omega\right) \left(\lambda - i\omega -\alpha_0 \right)}d \omega, & \\[0.2cm]
g_{22}(\lambda)=\dfrac{1}{2\pi}\dint_{-\infty}^{+\infty}
\dfrac{i\omega-a_{1}^{1}-b_{1}^{1}e^{-i\omega}}{h\left(i\omega\right)}\times & \\[0.2cm]
\qquad \qquad
\dfrac{-(a_{1}^{1}+a+\alpha_0)(\lambda-i\omega)-b + \alpha_0 a_{1}^{1}
-\left[(c+b_{1}^{1})(\lambda-i\omega)-\alpha_0 b_{1}^{1}+d\right]e^{-(\lambda-i\omega)}-r e^{-2(\lambda-i\omega)}}
{h\left(\lambda-i\omega\right) \left(\lambda - i\omega -\alpha_0 \right)}d \omega, & \\[0.2cm]
g_{31}(\lambda)=\dfrac{1}{2\pi}\dint_{-\infty}^{+\infty}
\dfrac{a_{2}^{1}+b_{2}^{1}e^{-i\omega}}{h\left(i\omega\right)}\times & \\[0.2cm]
\qquad \qquad
\dfrac{-(a_{2}^{2}+a+\alpha_0)(\lambda-i\omega)-b+\alpha_0 a_{2}^{2}
-\left[(c+b_{2}^{2})(\lambda-i\omega)-\alpha_0 b_{2}^{2}+d\right]e^{-(\lambda-i\omega)}-r e^{-2(\lambda-i\omega)}}
{h\left(\lambda-i\omega\right) \left(\lambda - i\omega -\alpha_0 \right)}d \omega, & \\[0.2cm]
g_{32}(\lambda)=\dfrac{1}{2\pi}\dint_{-\infty}^{+\infty}
\dfrac{a_{1}^{2}+b_{1}^{2}e^{-i\omega}}{h\left(i\omega\right)}\times & \\[0.2cm]
\qquad \qquad
\dfrac{-(a_{1}^{1}+a+\alpha_0)(\lambda-i\omega)-b + \alpha_0 a_{1}^{1}
-\left[(c+b_{1}^{1})(\lambda-i\omega)-\alpha_0 b_{1}^{1}+d\right]e^{-(\lambda-i\omega)}-r e^{-2(\lambda-i\omega)}}
{h\left(\lambda-i\omega\right) \left(\lambda - i\omega -\alpha_0 \right)}d \omega.
\end{array}\right.
}
\end{equation*}
Furthermore, there exist positive constants $c_{ij}\,(i=1,2,3,j=1,2)$ and $R_{ij}\,(i,j=1,2)$ such that for
$\Re(\lambda)>\alpha_0,$
\begin{equation}\label{Xij2}
\begin{array}{l}
\displaystyle\lim_{|\lambda|\to +\infty}|g_{ij}(\lambda)|=0\;(i=1,2,3,j=1,2), \\
\displaystyle \lim_{|\lambda|\to +\infty}\left|\lambda g_{12}(\lambda)\right|=c_{12},\quad
\lim_{|\lambda|\to +\infty}|h(\lambda-\alpha_0)g_{ii}(\lambda)|=c_{ii}\;(i=1,2), \\
\displaystyle \lim_{|\lambda|\to +\infty}\left|\lambda g_{21}(\lambda)\right|=c_{21},\quad
\lim_{|\lambda|\to +\infty}|h(\lambda-\alpha_0)g_{3j}(\lambda)|=c_{3j}\;(j=1,2) \\
\end{array}
\end{equation}
and
\begin{equation}\label{Xij3}
\left.\begin{array}{l}
\displaystyle
\lim_{|\lambda|\to +\infty}|\lambda\mathcal{L}((X_i^j)^2)(\lambda)|=R_{ij}\,(i,j=1,2),\\
\displaystyle \lim_{|\lambda|\to +\infty}|\lambda\mathcal{L}(X_1^1X_2^1)(\lambda)|=0,\quad
\lim_{|\lambda|\to +\infty}|\lambda\mathcal{L}(X_1^2X_2^2)(\lambda)|=0.
\end{array}\right.
\end{equation}
\end{lemma}
\Proof For any $\Re(\lambda)>2\alpha_{0}\,(\lambda\in\C),$
\begin{eqnarray*}
\mathcal{L}((X_1^1)^2)(\lambda)
&=&\dint_{0}^{+\infty}e^{-\lambda t}X_1^1(t)\dfrac{1}{2\pi}\int_{-\infty}^{+\infty}
\dfrac{i\omega-a_{2}^{2}-b_{2}^{2}e^{-i\omega}}{h(i\omega)}e^{i\omega t} d\omega dt \\
&=&\dfrac{1}{2\pi}\int_{-\infty}^{+\infty}\dfrac{i\omega-a_{2}^{2}-b_{2}^{2}e^{-i\omega}}{h(i\omega)}
\int_{0}^{+\infty}e^{-(\lambda -i\omega)t}X_1^1(t)dt d\omega\\
&=&\dfrac{1}{2\pi}\int_{-\infty}^{+\infty}
\dfrac{i\omega-a_{2}^{2}-b_{2}^{2}e^{-i\omega}}{h(i\omega)}
\dfrac{\lambda-i\omega-a_{2}^{2}-b_{2}^{2}e^{-(\lambda -i\omega)}}{h(\lambda-i\omega)}d\omega.
\end{eqnarray*}
Since $ \int_{0}^{+\infty}e^{-(\lambda - i\omega -\alpha_0)t}dt=\dfrac{1}{\lambda - i\omega -\alpha_0}$
for $\Re(\lambda)> \alpha_0$, then
\begin{eqnarray*}
\mathcal{L}((X_1^1)^2)(\lambda)
&=&\dfrac{1}{2\pi}\int_{-\infty}^{+\infty}
\dfrac{i\omega-a_{2}^{2}-b_{2}^{2}e^{-i\omega}}{h(i\omega)}\dfrac{1}{\lambda - i\omega -\alpha_0}d\omega\\
&&+\dfrac{1}{2\pi}\int_{-\infty}^{+\infty}\dfrac{i\omega-a_{2}^{2}-b_{2}^{2}e^{-i\omega}}{h(i\omega)}
\left[\dfrac{\lambda-i\omega-a_{2}^{2}-b_{2}^{2}e^{-(\lambda -i\omega)}}{h(\lambda-i\omega)}-\dfrac{1}{\lambda-i\omega-\alpha_0}\right]d\omega\\
&=&\int_{0}^{+\infty}e^{-(\lambda-\alpha_0)t} \dfrac{1}{2\pi}\int_{-\infty}^{+\infty}
\dfrac{i\omega-a_{2}^{2}-b_{2}^{2}e^{-i\omega}}{h(i\omega)}e^{i\omega t}d\omega dt+g_{11}(\lambda) \\
&=& \int_{0}^{+\infty}e^{-(\lambda-\alpha_0) t}X_1^1(t)dt+ g_{11}(\lambda)\\
&=& \dfrac{\lambda-\alpha_0-a_{}^{2}-b_{2}^{2}e^{-(\lambda-\alpha_0)}}{h(\lambda-\alpha_0)} + g_{11}(\lambda)
\end{eqnarray*}
for $\Re(\lambda) > \alpha_0$, where
\begin{equation*}
\footnotesize{
\left.\begin{array}{ll}
\displaystyle g_{11}(\lambda)=\dfrac{1}{2\pi}\int_{-\infty}^{+\infty}
\dfrac{i\omega-a_{2}^{2}-b_{2}^{2}e^{-i\omega}}{h\left(i\omega\right)}\times & \\[0.2cm]
\qquad\qquad \dfrac{-(a_{2}^{2}+a+\alpha_0)(\lambda-i\omega)-b+\alpha_0 a_{2}^{2}
-\left[(c+b_{2}^{2})(\lambda-i\omega)-\alpha_0 b_{2}^{2}+d\right]e^{-(\lambda-i\omega)}-r e^{-2(\lambda-i\omega)}}
{h\left(\lambda-i\omega\right)\left(\lambda - i\omega -\alpha_0 \right)}d \omega.
\end{array}\right.
}
\end{equation*}
Note that $g_{11}(\lambda)$ is convergent for any $\Re(\lambda)>\alpha_0.$
Hence there exist two positive constants $c_{11}$ and $R_{11}$ such that when $\Re(\lambda)>\alpha_0$,
\begin{equation*}
\lim_{|\lambda|\to +\infty}|g_{11}(\lambda)|=0, \;
\lim_{|\lambda|\to +\infty}|h(\lambda-\alpha_0)g_{11}(\lambda)|=c_{11},\;
\lim_{|\lambda|\to +\infty}|\lambda\mathcal{L}((X_1^1)^2)(\lambda)|=R_{11}.
\end{equation*}
Other expressions in \eqref{Xij1-1}-\eqref{Xij1-2}, \eqref{Xij2} and \eqref{Xij3} can be obtained similarly.
\qed

Similar to Lemma \ref{LemLyij}, we have the following expressions and estimates.
\begin{lemma}\label{LemLyij2}
Assume $\alpha_0 <0$. For any $ \Re(\lambda)>\alpha_0\,(\lambda\in\C),$
\footnotesize {
\begin{eqnarray*}\label{Xij4-1}
\mathcal{L}(X_1^1(t)X_1^1(t-1))(\lambda)&=&\dfrac{1}{2\pi}\int_{-\infty}^{+\infty}
\dfrac{e^{-i\omega}\left(i\omega-a_{2}^{2}-b_{2}^{2}e^{-i\omega}\right)}{h(i\omega)}
\cdot\dfrac{\lambda-i\omega-a_{2}^{2}-b_{2}^{2}e^{-(\lambda-i\omega)}}{h(\lambda-i\omega)} d\omega  \\
&=&\dfrac{e^{-(\lambda-\alpha_0)}\left(\lambda-\alpha_0-a_{2}^{2}-b_{2}^{2}e^{-(\lambda-\alpha_0)}\right)
}{h(\lambda-\alpha_0)}+\widetilde{g}_{11}(\lambda),\\
\mathcal{L}(X_1^2(t)X_1^2(t-1))(\lambda)&=&\dfrac{1}{2\pi}\int_{-\infty}^{+\infty}
\dfrac{e^{-i\omega}\left(a_{1}^{2}+b_{1}^{2}e^{-i\omega}\right)}{h(i\omega)}
\cdot\dfrac{a_{1}^{2}+b_{1}^{2}e^{-(\lambda-i\omega)}}{h(\lambda-i\omega)} d\omega  \\
&=&\dfrac{e^{-(\lambda-\alpha_0)}(a_{1}^{2}+b_{1}^{2}e^{-(\lambda-\alpha_0)})}{h(\lambda-\alpha_0)}+\widetilde{g}_{12}(\lambda),\\
\mathcal{L}(X_2^1(t)X_2^1(t-1))(\lambda)&=&\dfrac{1}{2\pi}\int_{-\infty}^{+\infty}
\dfrac{e^{-i\omega}\left(a_{2}^{1}+b_{2}^{1}e^{-i\omega}\right)}{h(i\omega)}
\cdot\dfrac{a_{2}^{1}+b_{2}^{1}e^{-(\lambda-i\omega)}}{h(\lambda-i\omega)} d\omega  \\
&=&\dfrac{e^{-(\lambda-\alpha_0)}(a_{2}^{1}+b_{2}^{1}e^{-(\lambda-\alpha_0)})}{h(\lambda-\alpha_0)}+\widetilde{g}_{21}(\lambda),\\
\mathcal{L}(X_2^2(t)X_2^2(t-1))(\lambda)&=&\dfrac{1}{2\pi}\int_{-\infty}^{+\infty}
\dfrac{e^{-i\omega}\left(i\omega-a_{1}^{1}-b_{1}^{1}e^{-i\omega}\right)}{h(i\omega)}
\cdot\dfrac{\lambda-i\omega-a_{1}^{1}-b_{1}^{1}e^{-(\lambda-i\omega)}}{h(\lambda-i\omega)} d\omega \\
&=&\dfrac{e^{-(\lambda-\alpha_0)}\left(\lambda-\alpha_0-a_{1}^{1}-b_{1}^{1}e^{-(\lambda-\alpha_0)}\right)}{h(\lambda-\alpha_0)}+\widetilde{g}_{22}(\lambda),\\
\mathcal{L}(X_1^1(t)X_2^1(t-1))(\lambda)&=&\dfrac{1}{2\pi}\int_{-\infty}^{+\infty}
\dfrac{e^{-i\omega}\left(a_{2}^{1}+b_{2}^{1}e^{-i\omega}\right)}{h(i\omega)}
\cdot\dfrac{\lambda-i\omega-a_{2}^{2}-b_{2}^{2}e^{-(\lambda-i\omega)}}{h(\lambda-i\omega)} d\omega  \\
&=&\dfrac{e^{-(\lambda-\alpha_0)}(a_{2}^{1}+b_{2}^{1}e^{-(\lambda-\alpha_0)})}{h(\lambda-\alpha_0)}+\widetilde{g}_{31}(\lambda),\\
\mathcal{L}(X_2^1(t)X_1^1(t-1))(\lambda)&=&\dfrac{1}{2\pi}\int_{-\infty}^{\infty}
\dfrac{e^{-i\omega}\left(i\omega-a_{2}^{2}-b_{2}^{2}e^{-i\omega}\right)}{h(i\omega)}
\cdot\dfrac{a_{2}^{1}+b_{2}^{1}e^{-(\lambda-i\omega)}}{h(\lambda-i\omega)}d\omega  \\
&=&\dfrac{e^{-(\lambda-\alpha_0)}\left(\lambda-\alpha_0-a_{2}^{2}-b_{2}^{2}e^{-(\lambda-\alpha_0)}\right)
}{h(\lambda-\alpha_0)}+\widetilde{g}_{32}(\lambda),\\
\mathcal{L}(X_1^2(t)X_2^2(t-1))(\lambda)&=&\dfrac{1}{2\pi}\int_{-\infty}^{+\infty}
\dfrac{e^{-i\omega}\left(i\omega-a_{1}^{1}-b_{1}^{1}e^{-i\omega}\right)}{h(i\omega)}
\cdot\dfrac{a_{1}^{2}+b_{1}^{2}e^{-(\lambda-i\omega)}}{h(\lambda-i\omega)} d\omega  \\
&=&\dfrac{e^{-(\lambda-\alpha_0)}\left(\lambda-\alpha_0-a_{1}^{1}-b_{1}^{1}e^{-(\lambda-\alpha_0)}\right)
}{h(\lambda-\alpha_0)}+\widetilde{g}_{41}(\lambda),\\
\mathcal{L}(X_2^2(t)X_1^2(t-1))(\lambda)&=&\dfrac{1}{2\pi}\int_{-\infty}^{+\infty}
\dfrac{e^{-i\omega}\left(a_{1}^{2}+b_{1}^{2}e^{-i\omega}\right)}{h(i\omega)}
\cdot\dfrac{\lambda-i\omega-a_{1}^{1}-b_{1}^{1}e^{-(\lambda-i\omega)}}{h(\lambda-i\omega)}d\omega  \\
&=&\dfrac{e^{-(\lambda-\alpha_0)}(a_{1}^{2}+b_{1}^{2}e^{-(\lambda-\alpha_0)})}{h(\lambda-\alpha_0)}+\widetilde{g}_{42}(\lambda),
\end{eqnarray*}
}
where
\begin{equation*}
\scriptsize{
\left.\begin{array}{ll}
\widetilde{g}_{11}(\lambda)=\dfrac{1}{2\pi}\dint_{-\infty}^{+\infty}
\dfrac{e^{-i\omega}\left(i\omega-a_{2}^{2}-b_{2}^{2}e^{-i\omega}\right)}{h\left(i\omega\right)}\times &\\[0.2cm]
\qquad\qquad
\dfrac{-(a_{2}^{2}+a+\alpha_0)(\lambda-i\omega)-b+\alpha_0 a_{2}^{2}
-\left[(c+b_{2}^{2})(\lambda-i\omega)-\alpha_0 b_{2}^{2}+d\right]e^{-(\lambda-i\omega)}-r e^{-2(\lambda-i\omega)}}
{h\left(\lambda-i\omega\right) \left(\lambda - i\omega -\alpha_0 \right)}d \omega,&\\[0.2cm]
\widetilde{g}_{12}(\lambda)=\dfrac{1}{2\pi}\dint_{-\infty}^{+\infty}
\dfrac{e^{-i\omega}\left(a_{1}^{2}+b_{1}^{2}e^{-i\omega}\right)}{h\left(i\omega\right)}\times  &\\[0.2cm]
\qquad\qquad
\dfrac{-(\lambda-i\omega)^2+(a_{1}^{2}-a)(\lambda-i\omega)-b-\alpha_0 a_{1}^{2}
-\left[(c-b_{1}^{2})(\lambda-i\omega)+\alpha_0 b_{1}^{2}+d\right]e^{-(\lambda-i\omega)}-r e^{-2(\lambda-i\omega)}}
{h\left(\lambda-i\omega\right) \left(\lambda - i\omega -\alpha_0 \right)}d \omega,&\\[0.2cm]
\widetilde{g}_{21}(\lambda)=\dfrac{1}{2\pi}\dint_{-\infty}^{+\infty}
\dfrac{e^{-i\omega}\left(a_{2}^{1}+b_{2}^{1}e^{-i\omega}\right)}{h\left(i\omega\right)}\times  &\\[0.2cm]
\qquad\qquad
\dfrac{-(\lambda-i\omega)^2+(a_{2}^{1}-a)(\lambda-i\omega)-b-\alpha_0 a_{2}^{1}
-\left[(c-b_{2}^{1})(\lambda-i\omega)+\alpha_0 b_{2}^{1}+d\right]e^{-(\lambda-i\omega)}-r e^{-2(\lambda-i\omega)}}
{h\left(\lambda-i\omega\right) \left(\lambda - i\omega -\alpha_0 \right)}d \omega, &\\[0.2cm]
\widetilde{g}_{22}(\lambda)=\dfrac{1}{2\pi}\dint_{-\infty}^{+\infty}
\dfrac{e^{-i\omega}\left(i\omega-a_{1}^{1}-b_{1}^{1}e^{-i\omega}\right)}{h\left(i\omega\right)}\times  &\\[0.2cm]
\qquad\qquad
\dfrac{-(a_{1}^{1}+a+\alpha_0)(\lambda-i\omega)-b + \alpha_0 a_{1}^{1}
-\left[(c+b_{1}^{1})(\lambda-i\omega)-\alpha_0 b_{1}^{1}+d\right]e^{-(\lambda-i\omega)}-r e^{-2(\lambda-i\omega)}}
{h\left(\lambda-i\omega\right) \left(\lambda - i\omega -\alpha_0 \right)}d \omega, &\\[0.2cm]
\widetilde{g}_{31}(\lambda)=\dfrac{1}{2\pi}\dint_{-\infty}^{+\infty}
\dfrac{e^{-i\omega}\left(a_{2}^{1}+b_{2}^{1}e^{-i\omega}\right)}{h\left(i\omega\right)}\times &\\[0.2cm]
\qquad\qquad
\dfrac{-(a_{2}^{2}+a+\alpha_0)(\lambda-i\omega)-b+\alpha_0 a_{2}^{2}
-\left[(c+b_{2}^{2})(\lambda-i\omega)-\alpha_0 b_{2}^{2}+d\right]e^{-(\lambda-i\omega)}-r e^{-2(\lambda-i\omega)}}
{h\left(\lambda-i\omega\right) \left(\lambda - i\omega -\alpha_0 \right)}d \omega, &\\[0.2cm]
\widetilde{g}_{32}(\lambda)=\dfrac{1}{2\pi}\dint_{-\infty}^{+\infty}
\dfrac{e^{-i\omega}\left(i\omega-a_{2}^{2}-b_{2}^{2}e^{-i\omega}\right)}{h\left(i\omega\right)}\times &\\[0.2cm]
\qquad\qquad
\dfrac{-(\lambda-i\omega)^2+(a_{2}^{1}-a)(\lambda-i\omega)-b-\alpha_0 a_{2}^{1}
-\left[(c-b_{2}^{1})(\lambda-i\omega)+\alpha_0 b_{2}^{1}+d\right]e^{-(\lambda-i\omega)}-r e^{-2(\lambda-i\omega)}}
{h\left(\lambda-i\omega\right) \left(\lambda - i\omega -\alpha_0 \right)}d \omega, &\\[0.2cm]
\widetilde{g}_{41}(\lambda)=\dfrac{1}{2\pi}\dint_{-\infty}^{+\infty}
\dfrac{e^{-i\omega}\left(i\omega-a_{1}^{1}-b_{1}^{1}e^{-i\omega}\right)}{h(i\omega)}\times &\\[0.2cm]
\qquad\qquad
\dfrac{-(\lambda-i\omega)^2+(a_{1}^{2}-a)(\lambda-i\omega)-b-\alpha_0 a_{1}^{2}
-\left[(c-b_{1}^{2})(\lambda-i\omega)+\alpha_0 b_{1}^{2}+d\right]e^{-(\lambda-i\omega)}-r e^{-2(\lambda-i\omega)}}
{h\left(\lambda-i\omega\right) \left(\lambda - i\omega -\alpha_0 \right)}d \omega, &\\[0.2cm]
\widetilde{g}_{42}(\lambda)=\dfrac{1}{2\pi}\dint_{-\infty}^{+\infty}
\dfrac{e^{-i\omega}\left(a_{1}^{2}+b_{1}^{2}e^{-i\omega}\right)}{h\left(i\omega\right)}\times &\\[0.2cm]
\qquad\qquad
\dfrac{-(a_{1}^{1}+a+\alpha_0)(\lambda-i\omega)-b + \alpha_0 a_{1}^{1}
-\left[(c+b_{1}^{1})(\lambda-i\omega)-\alpha_0 b_{1}^{1}+d\right]e^{-(\lambda-i\omega)}-r e^{-2(\lambda-i\omega)}}
{h\left(\lambda-i\omega\right) \left(\lambda - i\omega -\alpha_0 \right)}d \omega.
\end{array}\right.
}
\end{equation*}
Moreover, there exist positive constants $\widetilde{c}_{ij}\,(i=1,2,3,4,j=1,2)$ such that for $ \Re(\lambda)>\alpha_0,$
\begin{equation}\label{Xij5}
\lim_{|\lambda|\to +\infty}|\widetilde{g}_{ij}(\lambda)|=0\;(i=1,2,3,4, j=1,2),
\end{equation}
and when $i=1,3$
\begin{equation*}
\lim_{|\lambda|\to +\infty}|h(\lambda-\alpha_0)\widetilde{g}_{i1}(\lambda)|=\widetilde{c}_{i1},\;
\lim_{|\lambda|\to +\infty}\left|\lambda \widetilde{g}_{i2}(\lambda)\right|=\widetilde{c}_{i2},
\end{equation*}
when $i=2,4$,
\begin{equation}
\label{Xij5-2}
\lim_{|\lambda|\to +\infty}|h(\lambda-\alpha_0)\widetilde{g}_{i2}(\lambda)|=\widetilde{c}_{i2},\;
\lim_{|\lambda|\to +\infty}\left|\lambda \widetilde{g}_{i1}(\lambda)\right|=\widetilde{c}_{i1}.
\end{equation}
\end{lemma}

\section{Proof of Lemma \ref{Lem-Z}}
\Proof
Assume $\bar{\alpha}_0 <0$. Let
\begin{equation*}
Z(t) = \mathcal{L}^{-1}\left(\dfrac{1}{H(\lambda)}\mathrm{adj}(I-D(\lambda)) A(\lambda)\right)
\triangleq (Z_{ij}(t))_{2 \times 2}.
\end{equation*}
We only prove that there exists a positive constant $K_{11}$ such that for any $\beta>\max\{\bar{\alpha}_0, \beta_{0}\}$,
\begin{equation*}
|Z_{11}(t)| \leq K_{11}e^{\beta t},\; t>0.
\end{equation*}
The estimates of $Z_{12}(t),\;Z_{ii}(t)\; (i=1,2)$ are similar and omitted.

Let
$$
D(\lambda) = \left(\begin{array}{cc}
d_{11}(\lambda) \; & d_{12}(\lambda) \\[0.1cm]
d_{21}(\lambda) \; & d_{22}(\lambda)
\end{array}\right),
$$
then
\begin{equation} \label{Z11}
Z_{11}(t) = \mathcal{L}^{-1} \left(\dfrac{(1-d_{22})\mathcal{L}((X_1^1)^2)}{H(\lambda)}\right) +
\mathcal{L}^{-1} \left(\dfrac{d_{12}\mathcal{L}((X_2^1)^2)}{H(\lambda)}\right).
\end{equation}
Now, we estimate the two terms in \eqref{Z11} respectively.

When $\beta>\beta_{0}$, $H(\beta+i\omega)\neq 0$ for any $\omega\in\R$. Thus
\begin{eqnarray*}
\mathcal{L}^{-1}\left(\dfrac{(1-d_{22})\mathcal{L}((X_1^1)^2)}{H(\lambda)}\right)
&=&\int_{(c_{2})}\dfrac{(1-d_{22})\mathcal{L}((X_1^1)^2)}{H(\lambda)} e^{\lambda t} d\lambda\\
&=& \dfrac{1}{2\pi i}\lim_{T\to +\infty}\int_{c_{2}-iT}^{c_{2}+iT}
e^{\lambda t} \dfrac{(1-d_{22})\mathcal{L}((X_1^1)^2)}{H(\lambda)}d\lambda,
\end{eqnarray*}
where $c_{2}>\beta$ is large enough. First we want to prove that
\begin{eqnarray}
\mathcal{L}^{-1}\left(\dfrac{(1-d_{22})\mathcal{L}((X_1^1)^2)}{H(\lambda)}\right)
&=&\int_{(\beta)}e^{\lambda t}\dfrac{(1-d_{22})\mathcal{L}((X_1^1)^2)}{H(\lambda)}d\lambda.
\label{LW11}
\end{eqnarray}
To this end, we consider the integration of the function
$e^{\lambda t}\dfrac{(1-d_{22})\mathcal{L}((X_1^1)^2)}{H(\lambda)}$ around the boundary of the
box $\Gamma$ in the complex plan with the boundary $\Gamma_{1}\Gamma_{2}\Gamma_{3}\Gamma_{4}$ in the
anticlockwise direction, where the segment $\Gamma_{1}$ is the set $\{c_{2}+i\tau: -T\leq \tau \leq T\}$,
the segment $\Gamma_{2}$ is the set $\{u+iT: \beta\leq u \leq c_{2}\}$, the segment
$\Gamma_{3}$ is the set $\{\beta+i\tau: -T\leq \tau \leq T\}$ and the segment $\Gamma_{4}$ is the set
$\{u-iT: \beta\leq u \leq c_{2}\}$. Since $H(\lambda)$ has no zeros in this box $\Gamma$, the integral
over the boundary is zero, i.e.,
\begin{eqnarray*}
\left(\int_{\Gamma_{1}}+\int_{\Gamma_{2}}+\int_{\Gamma_{3}}+\int_{\Gamma_{4}}\right)
\left(e^{\lambda t}\dfrac{(1-d_{22})\mathcal{L}((X_1^1)^2)}{H(\lambda)}\right)d\lambda=0.
\end{eqnarray*}
Thus, \eqref{LW11} is concluded if
\begin{eqnarray*}
\lim_{T\to+\infty}\int_{\Gamma_{i}}e^{\lambda t}\dfrac{(1-d_{22})\mathcal{L}((X_1^1)^2)}{H(\lambda)}d\lambda=0\quad (i=2,4).
\end{eqnarray*}
Since $H(\lambda) =  (1-d_{11})(1-d_{22}) - d_{12}d_{21}$, we have
\begin{equation*}
\dfrac{H(\lambda)}{1-d_{22}} = 1-d_{11} -\dfrac{d_{12}d_{21}}{1-d_{22}}.
\end{equation*}
From Lemma \ref{Lem-D}, there exits a constant $T_{1}\geq T_{0}>0$ such that when
$|\lambda|\geq T_{1}$ and $\Re(\lambda)>\bar{\alpha}_0,$
\begin{equation*}
\left|\dfrac{H(\lambda)}{1-d_{22}}\right|
\geq  1-|d_{11}|-\dfrac{|d_{12}||d_{21}|}{1-|d_{22}|}
\geq 1-\dfrac{d_0}{|\lambda|}-\dfrac{ \dfrac{d_0^2}{|\lambda|^2} }{ 1-\dfrac{d_0}{|\lambda|} }
\geq  \dfrac{1}{2}.
\end{equation*}
Hence  from \eqref{Xij3}, for $|\lambda|\geq T_{1}$ and $\Re(\lambda)>\bar{\alpha}_0$,
$|\mathcal{L}((X_1^1)^2) \leq \dfrac{1+R_{11}}{|\lambda|}$, and therefore
\begin{eqnarray*}
\left|\int_{\Gamma_{2}}e^{\lambda t}\dfrac{(1-d_{22})\mathcal{L}((X_1^1)^2)}{H(\lambda)}d\lambda\right|
&=&\left|\int_{c_{2}+iT}^{\beta+iT}e^{\lambda t}\dfrac{(1-d_{22})\mathcal{L}((X_1^1)^2)}{H(\lambda)}d\lambda\right|\\[0.1cm]
&\leq & 2\int_{\beta}^{c_{2}}e^{u t}\dfrac{1+R_{11}}{\sqrt{u^2+T^2}}du \\
&\leq &  \dfrac{2e^{c_{2}t} (1+R_{11})}{T} (c_{2}-\beta) \to 0 \;(\text{as} \;T\to +\infty).
\end{eqnarray*}
Similarly,
\begin{equation*}
\int_{\Gamma_{4}}\dfrac{(1-d_{22})\mathcal{L}((X_1^1)^2)}
{H(\lambda)}e^{\lambda t}d\lambda\to 0 \;(\text{as} \;T\to +\infty).
\end{equation*}
Thus, \eqref{LW11} is obtained.

Let $T_{1}$ as above, and
\begin{eqnarray*}
W(\lambda)&=&\dfrac{(1-d_{22})\mathcal{L}((X_1^1)^2)}{H(\lambda)}-\dfrac{1}{\lambda-\beta_{0}}\\
&=&\dfrac{(1-d_{22})\mathcal{L}((X_1^1)^2)}{H(\lambda)}\cdot\dfrac{\lambda-\beta_{0}-\dfrac{H(\lambda)}{(1-d_{22})\mathcal{L}((X_1^1)^2)}}{\lambda-\beta_{0}}.
\end{eqnarray*}
From Lemma \ref{LemLyij}, we have for $\Re(\lambda)>\bar{\alpha}_0,$
\begin{eqnarray*}
&&\lambda-\beta_{0}-\dfrac{H(\lambda)}{(1-d_{22})\mathcal{L}((X_1^1)^2)}\\
&=& \lambda-\dfrac{1}{\mathcal{L}((X_1^1)^2)}-\beta_{0}
-\dfrac{1}{\mathcal{L}((X_1^1)^2)}\left(\dfrac{H(\lambda)}{ 1-d_{22} }-1\right) \\[0.1cm]
&=&-\beta_{0}-\dfrac{1}{\mathcal{L}((X_1^1)^2)}\left(\dfrac{H(\lambda)}{ 1-d_{22} }-1\right)
+\dfrac{(\alpha_0-a - a_2^2)\lambda +  \alpha_0 a -b -h(\lambda-\alpha_0)g_{11}}
{\lambda-\alpha_0-a_{2}^{2}-b_{2}^{2}e^{-(\lambda-\alpha_0)} + h(\lambda-\alpha_0) g_{11}} \\[0.1cm]
&&{}- \dfrac{ \left(b_2^2 \lambda +c(\lambda-\alpha_0)+ d \right)e^{-(\lambda-\alpha_0)}+ re^{-2(\lambda-\alpha_0)} } {\lambda-\alpha_0-a_{2}^{2}-b_{2}^{2}e^{-(\lambda-\alpha_0)} + h(\lambda-\alpha_0) g_{11}}
\end{eqnarray*}
and
\begin{eqnarray*}
\dfrac{1}{\mathcal{L}((X_1^1)^2)}\left(\dfrac{H(\lambda)}{ 1-d_{22} }-1\right)
&=& \dfrac{h(\lambda-\alpha_0)}{\lambda-\alpha_0-a_{2}^{2}-b_{2}^{2}e^{-(\lambda-\alpha_0)} + h(\lambda-\alpha_0) g_{11}}\times \\ &&{}\left(-d_{11}-\dfrac{d_{12}d_{21}}{1-d_{22}}\right).
\end{eqnarray*}
From \eqref{dij2}, there exists a positive constant $\widetilde{r}_{11}$ such that for
$|\lambda|\geq T_{1}$ and $\Re(\lambda)>\bar{\alpha}_0,$
\begin{eqnarray*}
\left|\dfrac{1}{\mathcal{L}((X_1^1)^2)}\left(\dfrac{H(\lambda)}{1-d_{22}}-1\right)\right|
&\leq& \widetilde{r}_{11}
\end{eqnarray*}
and
\begin{equation*}
\left|\lambda-\beta_{0}-\dfrac{H(\lambda)}{(1-d_{22})\mathcal{L}((X_1^1)^2)}\right|
\leq |\alpha_0-a - a_2^2|+ |\beta_{0}| + \widetilde{r}_{11}
= |\alpha_0+ a_1^1|+ |\beta_{0}| + \widetilde{r}_{11},
\end{equation*}
which, applying \eqref{Xij3}, implies
\begin{eqnarray*}
\left|W(\lambda)\right|&\leq& \dfrac{|\alpha_0+ a_1^1|+ |\beta_{0}| +\widetilde{r}_{11}}{|\lambda-\beta_{0}|} \cdot
\dfrac{2(1+R_{11})}{|\lambda|}
\leq \dfrac{2(1+R_{11})\left(|\alpha_0+ a_1^1|+ |\beta_{0}| +\widetilde{r}_{11}\right)}{|\lambda||\lambda-\beta_{0}|}.
\end{eqnarray*}
Thus, for $\beta>\max\{\bar{\alpha}_0, \beta_{0}\}$, since $\lambda-\beta_0\not=0$ and $H(\lambda)\not = 0\;(\Re(\lambda)>\beta_0)$,
\begin{eqnarray*}
\left|\mathcal{L}^{-1}\left(\dfrac{(1-d_{22})\mathcal{L}((X_1^1)^2)}{H(\lambda)}\right)\right|
&=&\left|\int_{(\beta)}e^{\lambda t}\dfrac{(1-d_{22})\mathcal{L}((X_1^1)^2)}{H(\lambda)}d\lambda\right| \\
&=&\left|\int_{(\beta)}e^{\lambda t} W(\lambda)d\lambda +\int_{(\beta)}\dfrac{e^{\lambda t}}{\lambda-\beta_{0}}d\lambda\right|\\
&\leq&e^{\beta_{0}t}+\left|\int_{(\beta)}e^{\lambda t} W(\lambda)d\lambda\right|\\
&\leq &  e^{\beta t}+ e^{\beta t} \Big(\int_{-T_1}^{T_1} |W(\beta + i\tau)| d\tau\\
&&{}\qquad + 2 \int_{T_1}^{+\infty} \dfrac{2(1+R_{11})\left(|\alpha_0+ a_1^1|+ |\beta_{0}|
+\widetilde{r}_{11}\right)}{\sqrt{\beta^2+\tau^2}\sqrt{(\beta-\beta_{0})^2+\tau^2}}d\tau \Big) \\
&\triangleq &\tilde{K}_{11}e^{\beta t}, \quad t>0,
\end{eqnarray*}
where
\begin{equation*}
\tilde{K}_{11}=1+\int_{-T_1}^{T_1} |W(\beta + i\tau)| d\tau + 2 \int_{T_1}^{+\infty}
\dfrac{2(1+R_{11})\left(|\alpha_0+ a_1^1|+ |\beta_{0}| +\widetilde{r}_{11}\right)}
{\sqrt{\beta^2+\tau^2}\sqrt{(\beta-\beta_{0})^2+\tau^2}}d\tau .
\end{equation*}

Now, we consider the second term in \eqref{Z11}. Since
\begin{equation*}
\dfrac{H(\lambda)}{d_{12}(\lambda)} = \dfrac{(1-d_{11})(1-d_{22})}{d_{12}}-d_{21},
\end{equation*}
from Lemma \ref{Lem-D}, there exits a constant $T_{2}\geq T_{0}>0$ such that when
$|\lambda|\geq T_{2}$ and $\Re(\lambda)>\bar{\alpha}_0,$
\begin{eqnarray*}
\left|\dfrac{H(\lambda)}{\lambda d_{12}(\lambda)}\right|
&\geq & \dfrac{1}{|\lambda|}\left[ \dfrac{(1-|d_{11}|)(1-|d_{22}|)}{|d_{12}|}-|d_{21}| \right] \\[0.1cm]
& \geq&   \dfrac{(1-\dfrac{d_0}{|\lambda|})(1-\dfrac{d_0}{|\lambda|})}{d_0}-\dfrac{d_0}{|\lambda|^2}
\geq \dfrac{1}{2d_0}.
\end{eqnarray*}
Thus for $|\lambda|\geq T_{2}, \,\Re(\lambda)>\bar{\alpha}_0$ and $i=2,4$, from \eqref{Xij3},
\begin{eqnarray*}
\left|\int_{\Gamma_{i}}e^{\lambda t}\dfrac{d_{12}\mathcal{L}((X_2^1)^2)}{H(\lambda)}d\lambda\right|
&\leq & \int_{\beta}^{c_{2}}e^{u t} \dfrac{ 2d_0 (1+R_{21})}{u^2 + T^2}du \\[0.1cm]
&\leq & \dfrac{2e^{c_{2}t} d_0 (1+R_{21}) }{T^2} (c_{2}-\beta) \to 0 \;(\text{as} \;T\to +\infty).
\end{eqnarray*}
Hence
\begin{eqnarray*}
\mathcal{L}^{-1}\left(\dfrac{d_{12}\mathcal{L}((X_2^1)^2)}{H(\lambda)}\right)
&=&\int_{(\beta)}e^{\lambda t}\dfrac{d_{12}\mathcal{L}((X_2^1)^2)}{H(\lambda)}d\lambda.
\end{eqnarray*}
Therefore, from \eqref{Xij3}, we obtain for $\beta>\max\{\bar{\alpha}_0, \beta_{0}\}$,
\begin{eqnarray*}
\left|\mathcal{L}^{-1}\left(\dfrac{d_{12}\mathcal{L}((X_2^1)^2)}{H(\lambda)}\right)\right|
& \leq & e^{\beta t} \left(\int_{-T_2}^{T_2} \left|\dfrac{d_{12}\mathcal{L}((X_2^1)^2)}{H(\beta + i\tau)}\right| d\tau
+ 2 \int_{T_2}^{+\infty}\dfrac{2d_0 (1+R_{21})}{\beta^{2} + \tau^{2}}d\tau \right) \\
& \triangleq & \bar{K}_{11}e^{\beta t}, \quad t>0,
\end{eqnarray*}
where
\begin{equation*}
\bar{K}_{11}=\int_{-T_2}^{T_2} \left|\dfrac{d_{12}\mathcal{L}((X_2^1)^2)}{H(\beta + i\tau)}\right| d\tau
+ 2 \int_{T_2}^{+\infty}\dfrac{2d_0 (1+R_{21})}{\beta^{2} + \tau^{2}}d\tau .
\end{equation*}
Taking $K_{11}=\tilde{K}_{11} + \bar{K}_{11}$, \eqref{Z11} is concluded and the Lemma is proved.
\qed

\end{appendix}

\section*{References}
\bibliographystyle{elsarticle-num}

\begin{thebibliography}{99}

\bibitem{Arino} O. Arino, M. L. Hbid and E. Ait Dads, { Delay Differential Equations
and Application}, Springer, 2006.


\bibitem{Bellman} R. Bellman and K. L. Cooke, { Differential-Difference Equations},
Academic, New York, 1963.


\bibitem{Caraballo10} T. Caraballo, J. Duan, K. Lu and B. Schmalfu{\ss}, { Invariant manifolds
for random and stochastic partial differential equations}, Advanced Nonlinear Studies, {\bf 10}
(2010) 23-52.


\bibitem{Duan04} J. Duan, K. Lu, and B. Schmalfu{\ss}, { Smooth stable and unstable manifolds for
stochastic evolutionary equations}, J. Dynam. and Diff. Eqns., {\bf 16} (2004) 949-972.


\bibitem{Hale93} J. K. Hale and S. M. Verduyn Lunel, {  Introduction to Functional
Differential Equations}, Springer Press, New York, 1993.


\bibitem{Ito} K. It\^{o} and M. Nisio, { On stationary solutions of a stochastic
differential equations}, J. Math. Kyoto Univ., {\bf 4} (1964) 1-75.


\bibitem{Ivanov03} A. F. Ivanov, Y. I. Kazmerchuk and A. V. Swishchuk, { Theory,
stochastic stability and applications of stochastic delay differential equations:
A survey of results}, Differential Equations Dynamical Systems, {\bf 11} (2003) 55-115.


\bibitem{LeiM} J. Lei and M. C. Mackey, { Stochastic differential delay equation,
Moment stability and its application to the hamatopoietic stem cell regulation system},
SIAM J. Appl. Math., {\bf 67} (2007) 387-407.

\bibitem{Mackey94} M. C. Mackey and I. G. Nechaeva, { Noise and Stability in Differential Delay Equations},
J. Dynam. and Diff. Eqns., {\bf 6} (1994) 395-426.


\bibitem{Mackey95} M. C. Mackey and I. G. Nechaeva, { Solution moment stability in stochastic differential
delay equations}, Phsical Review E, {\bf 52} (1995) 3366-3376.


\bibitem{Mao1} X. Mao, { Stochastic Differential Equations and Their Applications},
Horwood Publishing, Chichester, UK, 1997.

\bibitem{Mao2} X. Mao and S. Sabanis, { Numerical solutions of stochastic differential
delay equations under local Lipschitz condition}, J. Comput. Appl. Math., {\bf 151} (2003) 215-227.


\bibitem{Mao3} X. Mao, { Attraction, stability and boundedness for stochastic
differential delay equations}, Nonlinear Analysis,  {\bf 47} (2001) 4795-4806.


\bibitem{Mohammed84} S.-E. A. Mohammed, { Stochastic Functional Differential Equations},
Res. Notes in Math. 99, Pitman, Boston, 1984.


\bibitem{Mohammed03} S.-E. A. Mohammed and M. K. R. Scheutzow, { The stable manifold theorem for non-linear
stochastic systems with memory. I. Existence of the semiflow}, Journal of Functional Analysis,
 {\bf 205} (2003) 271-305.


\bibitem{Mohammed04} S.-E. A. Mohammed and M. K. R. Scheutzow, { The stable manifold theorem for non-linear
stochastic systems with memory II. The local stable manifold theorem}, Journal of Functional Analysis,
 {\bf 206} (2004) 253-306.


\bibitem{WLL} Z. Wang, X. Li and J. Lei, { Moment boundedness of linear stochastic differential equation with distributed delay},
\textit{Preprint submitted to Stochastic Processes and their Applications}.


\end{thebibliography}

\end{document}